\documentclass[10pt,a4paper,oneside,leqno]{article}
\usepackage{amsfonts,amssymb}
\usepackage{amscd}
\newtheorem{defn}{Definition}

\newtheorem{lemma}{Lemma}

\newtheorem{observen}{Observation}

\vspace{0.5cm}
 4
\font\ebf=cmbx8
\font\erm=cmr8

\usepackage{graphicx, graphics}
\usepackage{color}
\newcommand{\os}{\oplus\!\!\to}

\newcommand{\layer}[2]{\langle\Phi_#1 \to \Phi_#2\rangle}
\newcommand{\fnomial}[2]{ {{#1} \choose {#2}}_F }

\begin{document}
\begin{center}
	\noindent { \textsc{ On natural join of posets properties and first applications.}}  \\ 
	\vspace{0.3cm}
	\vspace{0.3cm}
	\noindent Andrzej Krzysztof Kwa\'sniewski \\
	\vspace{0.2cm}
	\noindent {\erm Member of the Institute of Combinatorics and its Applications  }\\
{\erm High School of Mathematics and Applied Informatics} \\
	{\erm  Kamienna 17, PL-15-021 Bia\l ystok, Poland }\noindent\\
	\noindent {\erm e-mail: kwandr@gmail.com}\\
\vspace{0.2cm}
\end{center}

\noindent {\ebf Abstract:}
\vspace{0.1cm}
\noindent {\small In early beginnings of the past century Felix Hausdorff introduced the concept of a partially ordered set thus extending Richard Dedekind lattice theory which began in the early 1890s.  Then the subject lay more or less  dormant until Garrett Birkhoff, Oystein Ore and others picked it up in the 1930s. Since then, many noted mathematicians have contributed to the subject, including Garrett Birkhoff, Richard Dedekind, Israel Gelfand, George Grätzer, Aleksandr Kurosh, Anatoly Malcev, Oystein Ore, Gian-Carlo Rota, Alfred Tarski and Johnny von Neumann  (see Garrett Birkhoff [1] for necessary relevant definitions - including operations on posets). Now poset and lattice theories are powerful theories of accelerated  growth inseparably tied  with graph theory. \\ 
Here - in addition to the three operations on posets from [1] - one introduces the natural join $P \os Q$ of posets $\left\langle P,\leq_P \right\rangle$ and $\left\langle Q,\leq_Q \right\rangle$ which should satisfy the  following \textcolor{blue}{\textbf{Natural Join Conformity Condition}}:  partial orders $\leq_P $  and $\leq_Q$ must be  equivalent on corresponding nonempty   $R \subseteq P$  and $R'\subseteq Q$ isomorphic sub-posets $R$ and $R'$ which we then shall consider identical by convention. In particular cases  $P \os Q$ may be expressed via ordinal sum  $\oplus $ accompanied with projection out operation as it is the case with sequences of antichains forming all together $r$-ary complete ($\equiv$ universal) relations via cardinal sum of these trivial posets. For other, earlier definitions of natural join of graded posets in terms of natural join  of bipartite natural join summands of graded poset $P$ or via natural join of corresponding  bi-adjacency matrices of bipartite natural join of these summands - see [2,3,4,5,6] or  for example  Definitions 12 or 2  in this paper.\\
The present definition of natural join $P \os Q$ of posets offers many conveniences. For example one arrives at a very simple proof of the M{\"{o}}bius function formula for cobweb posets.\\
The main aims of this article are to present at first the authors update applications of natural join of posets and/or  their cover relation Hasse digraphs and then one summarizes more systematically the general properties of natural join of posets with posing some questions arising on the way.\\
Thus in this note - apart from revealing some general properties of natural join of posets we also quote the authors combinatorial interpretations of cobweb posets and  explicit formulas  for the  zeta matrix [2-23] as well as the  inverse  of zeta matrix for \textbf{any} graded posets [2,3] using Knuth notation [24] while following [3,2].\\
These formulas use the supplemented by sub-matrix factors  corresponding  formulas for cobweb posets and their Hasse digraphs named KoDAGs, which are interpreted also as chains of binary \textbf{complete} (or universal) relations -  joined by the natural join operation  - see [2,3,4,5,6,12].\\
Such cobweb posets and equivalently their Hasse digraphs named KoDAGs -  are also encoded by discrete hyper-boxes [12] and the natural join operation of such discrete hyper-boxes is just Cartesian product of them accompanied with projection out of - sine qua non - common faces [2]. All graded posets with no mute vertices in their Hasse diagrams [2] (i.e. no "internal" vertex has  in-degree or out-degree equal zero - see Fig.2) are natural join of chain of relations and may be at the same time  interpreted as an $n-ary$, $n \in N \cup \left\{\infty \right\}$ relation. The Whitney numbers and characteristic polynomials explicit formulas for  cobweb posets are supplied too, following their derivation in [2]. Further reading on cobweb posets, KoDags and related stuff - see [25-49]. Knuth notation [24] and the author's upside down notation  are used throughout the paper.
}

\vspace{0.2cm}

\noindent Key Words: graded digraphs, cobweb  posets,  natural join, ordinal sum 

\vspace{0.1cm}

\noindent AMS Classification Numbers: 06A06 ,05B20, 05C75  

\vspace{0.1cm}

\noindent  This is The Internet Gian-Carlo Rota Polish Seminar article, No \textbf{9}, \textbf{Subject 4}, 2009 August 10\\
\emph{http://ii.uwb.edu.pl/akk/sem/sem\_rota.htm}\\

%%%%%%%%%%%%%%%%%%%%%%%%%%%%%%%%%%%%%%%%%%%%%%%%%%%%%%%%%%%%%%%%%%
%%%%%%%%%%%%%%%%%%%%%%%%%%%%%%%%%%%%%%%%%%%%%%%%%%%%%%%%%%%%%%%%%%

\section{Properties of natural join and application to cobwebs  $\Pi$ including M{\"{o}}bius  function  formula.} 

%%%%%%%%%%%%%%%%%%%%%%%%%%%%%%%%%%%%%%%%%%%%%%%%%%%%%%%%%%%%%%%%%%%%

\noindent \textbf{1.1.} \textbf{Ponderables.} [6,5,4,3,2]\\
%%%%%%%%%%%%%%%%%%%%%%%%%%%%%%%%%%%%%%%%%%%%%%%%%%%%%%%%%%%%%%%%%%
%%%%%%%%%%%%%%%%%%%%%%%%%%%%%%%%%%%%%%%%%%%%%%%%%%%%%%%%%%%%%%%%%%

\vspace{0.1cm}
\noindent We shall here take for granted the notation and the results of [5,4,3,2]. In particular $\left\langle \Pi,\leq \right\rangle$ denotes cobweb partial order set (cobweb poset)  while   $I(\Pi,R)$ denotes its incidence algebra over the ring $R$.  Correspondingly  $\left\langle P,\leq \right\rangle$ denotes arbitrary graded poset while   $I(P,R)$ denotes its incidence algebra over the ring $R$. 
$[n]\equiv \left\{1,2,...,n\right\}$ , for example $[k_F]\equiv \left\{1,2,...,k_F\right\}$.
\vspace{0.1cm}
\noindent $R$ might be taken to be Boolean algebra $2^{\left\{1 \right\}}$ , the field $Z_2=\left\{0,1\right\}$, the ring of integers $Z$ or real or complex or $p$-adic fields. The present article is the next one in a series of papers listed in  order of appearence and these are:[6], [5],[4],[3],[2]. The abbreviation \textbf{DAG} $\equiv$ \textbf{D}irected \textbf{A}cyclic \textbf{G}raph. \\
\vspace{0.1cm}
\noindent Inspired by Gaussian integers sequence notation $\left\{ n_q \right\}_{n\geq 0}$ - the authors   upside down notation  is used throughout this paper i.e. $F_n \equiv  n_F$. 
\textcolor{blue}{\textbf{The Upside Down Notation }} was used since last century effectively (see [2-22] and [29-48]).\\
Through all the paper  $F$ denotes a natural numbers valued sequence $\left\{ n_F \right\}_{n\geq 0} \equiv  \left\{ F_n \right\}_{n\geq 0}$  sometimes specified to be Fibonacci or others - if needed. Among many consequences of this is that \textbf{\textbf{graded posets}} ($\equiv$ their cover relation  digraphs $ \Longleftrightarrow$ Hasse diagrams) \textbf{\textit{are connected}} and sets of their minimal elements are finite.\\
\vspace{0.1cm}
\begin{defn}
\noindent Let  $F = \left\langle k_F \right\rangle_{k=0}^n$ be an arbitrary natural numbers valued sequence, where $n\in N \cup \left\{0\right\}\cup \left\{\infty\right\}$. We say that the graded poset $P = (\Phi,\leq)$, where $P= \bigcup_{k=0}^n \Phi_k$ is \textcolor{red}{\textbf{denominated}} (encoded=labeled) by  $F$  iff   $\left|\Phi_k\right| = k_F$ for $k = 0,1,..., n.$ . We shall also use the expression - $F$-graded poset.
\end{defn}

\begin{defn}
\noindent  Let  $n\in N \cup \left\{0\right\}\cup \left\{\infty\right\}$. Let   $r,s \in N \cup \left\{0\right\}$.  Let  $\Pi_n$ be the graded partial ordered set (poset) i.e. $\Pi_n = (\Phi_n,\leq)= ( \bigcup_{k=0}^n \Phi_k ,\leq)$ and $\left\langle \Phi_k \right\rangle_{k=0}^n$ constitutes ordered partition of $\Pi_n$. A graded poset   $\Pi_n$  with finite set of minimal elements is called \textbf{cobweb poset} \textsl{iff}  
$$\forall x,y \in \Phi \  i.e. \  x \in \Phi_r \ and \  y \in \Phi_s \ :  r \neq s\ \Rightarrow \   x < y   \ or \ y < x  , $$ 
$\Pi_\infty \equiv \Pi. $
\end{defn}
\noindent \textbf{Note}. By definition of $\Pi_n$ being graded its  levels $\Phi_r \in \left\{\Phi_k\right\}_k^n$ are independent sets  $\equiv$  antichains, $n\in N \cup \left\{0\right\}\cup \left\{\infty\right\}$. 
\vspace{0.2cm}
\noindent The Definition 2  is the reason for calling Hasse digraph $D = \left\langle \Phi, \leq \cdot \right\rangle $ of the poset $\Pi = (\Phi,\leq))$ a \textbf{\textcolor{red}{Ko}}DAG as in  Professor \textbf{\textcolor{red}{K}}azimierz   \textbf{\textcolor{red}{K}}uratowski native language one word \textbf{\textcolor{red}{Ko}mplet} means \textbf{complete ensemble}- see more in  [3,4] and for the history of this name see:  The Internet Gian-Carlo Polish Seminar Subject 1. \textit{ oDAGs and KoDAGs in Company }(Dec. 2008). Examples - see Fig.1, Fig.3  and Fig.11; note substantialy different Fig.2 and consult also  [2-6] for other figures and references therein.

\vspace{0.1cm}
%%%%%%%%%%%%%%%%%%%%%%%%%%%%%%%%%%%%%%%%%%%%%%%%%%%%%%%%%%%%%%%%%%

\vspace{0.2cm}
\begin{figure}[ht]
\begin{center}
	\includegraphics[width=100mm]{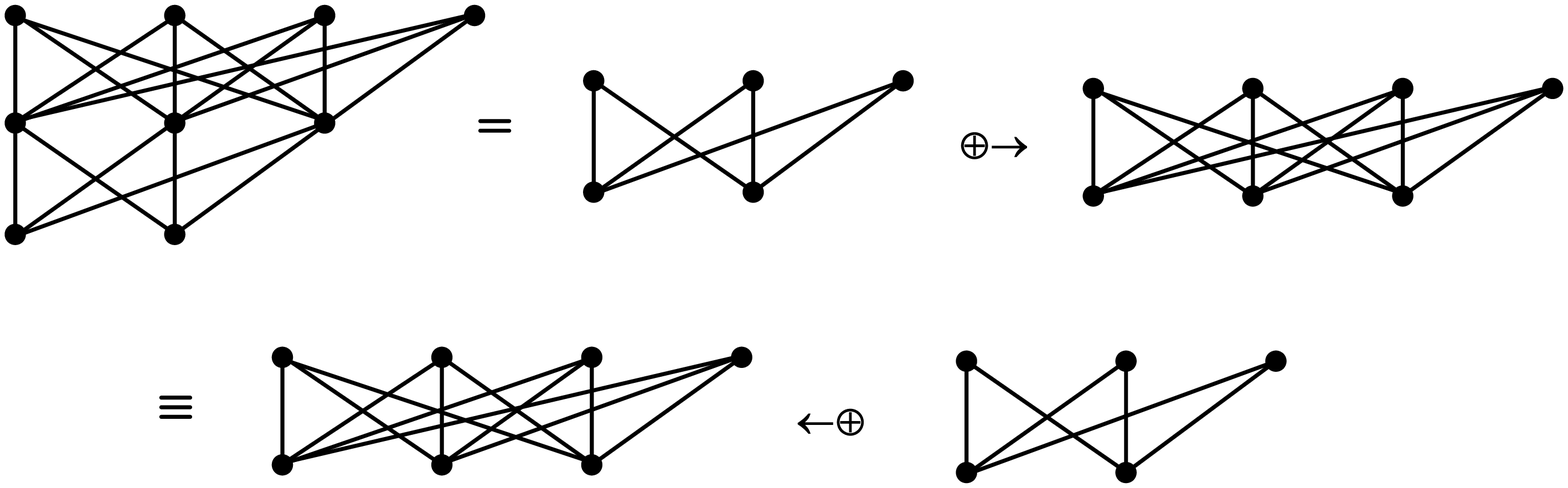}
	\caption 
	{Display of the natural join of bipartite layers $\left\langle \Phi_k \rightarrow \Phi_{k+1} \right\rangle$ $F = N $.} \label{fig:representation}
\end{center}
\end{figure}

\vspace{0.1cm}

\noindent \textbf{Note}.  The graded digraph   $G$  (graph)  is utterly denominated by its sequence of  bipartite digraphs (graphs)  that every two consecutive  levels of  $G$ do constitute.\\ 
The complete graded digraph (KoDAG) denoted as $K$  is utterly denominated by its sequence of complete bipartite digraphs that every two consecutive  levels of  $K$  do constitute.\\ 
Because of their appearance, an at a first glance outlook - posets characterized by these complete  graded  digraphs  $K$ where called cobweb posets i.e. digraphs $K$ are identified with Hasse diagrams of cobweb posets.\\
This appearance of complete or almost complete graded digraphs is tremendously prevailing - like the hoary trees with silver cobweb threads.  The opposite to this picture of a tree with cobweb  is a bare rooted  tree graph - void of this spider's web  hoary tunicate; see Fig.3 and figures in [2-15], [\textbf{17}], [18-23].

\vspace{0.1cm}

\noindent \textbf{Note [4]}. Cobweb posets graphs are examples of  Differential posets  due to Stanley  [50] and might serve for  more general structures  [due to Fomin [51] ] -  called  Dual graded Graphs. This is accompanied with a simple though fascinating observation [4] :  structure defining operators associated  to Dual graded Graphs generate  GHW algebra! (GHW = Graves-Heisemberg-Weil; see [44,33,30,31,52]).

\vspace{0.1cm}	

\begin{figure}[ht]
\begin{center}
	\includegraphics[width=100mm]{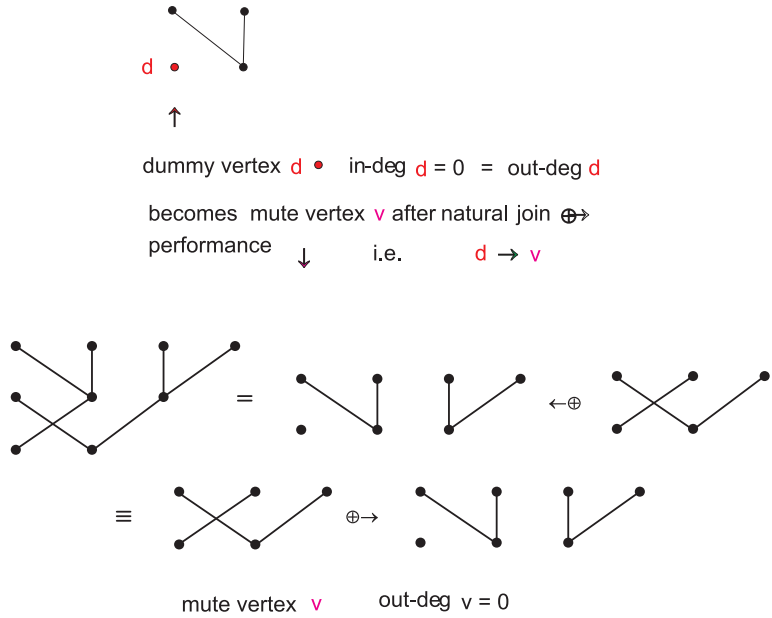}
	\caption {Display of the  natural join $\os$ of  bipartite digraphs with one dummy vertex resulting in one mute node. }\end{center}\label{fig:representation}
\end{figure}

\vspace{0.1cm}

\vspace{0.2cm}

\noindent \textbf{1.2.} \textbf{Ordinal and cardinal sums with added natural join of posets.} 

%%%%%%%%%%%%%%%%%%%%%%%%%%%%%%%%%%%%%%%%%%%%%%%%%%%%%%%%%%%%%%%%%%
%%%%%%%%%%%%%%%%%%%%%%%%%%%%%%%%%%%%%%%%%%%%%%%%%%%%%%%%%%%%%%%%%%
\vspace{0.2cm}
\noindent  Let us recall that the ordinal sum [linear sum] of two disjoint ordered sets $P$ and $Q$, denoted by  $P\oplus Q$, is the union of $P$ and $Q$, with $P$'s elements ordered as in $P$ while $Q$'s elements are correspondingly ordered as in $Q$, and  for each $x\in P$ and $y \in Q$ we put  $x \leq y$ .  
\noindent  The Hasse diagram of  $P\oplus Q$ we construct placing $Q$'s diagram just above $P$'s diagram and with an edge between each minimal element of $Q$ and each maximal element of $P$. 
(Definitions and properties of ordinal sum etc. see: [1] Chapter  III). 
\vspace{0.2cm}

\noindent Here we propose to add to  the standard operations   $1,2,3$  on ordered sets from [1]  i.e. to  operations  below \\
\vspace{0.1cm}
\noindent\textbf{ 1.2.1.}     dual  $P^*$  of  $P$ ;\\
\vspace{0.1cm}
\noindent \textbf{1.2.2.}     the disjoint union  $P+Q =$  cardinal sum ;\\          
\vspace{0.1cm}
\noindent \textbf{1.2.3.}     the ordinal sum   $P\oplus Q$ .\\
\vspace{0.1cm} 
\noindent  - the fourth one. Namely - we propose to \textbf{add} the following next binary operation on  posets - especially comfortable and handy  for \textbf{graded} posets.\\
\vspace{0.1cm}
\noindent \textbf{1.2.4.} This fourth one we propose as in [5,4,3,2] is the  \textcolor{red}{\textbf{natural join}} $P \os Q$, here below defined via \textbf{5.d}  which for particular cases is expressed by ordinal sums  $P\oplus Q$  and/or disjoint sums.

\vspace{0.1cm}

\noindent In order to define  natural join $P \os Q$ of posets $\left\langle P,\leq_P \right\rangle$ and $\left\langle Q,\leq_Q \right\rangle$  both arguments of the $\os$ operation should satisfy the \textcolor{blue}{\textbf{Natural Join Conformity Condition}} which is the straightforward generalization of the Domain-Codomain $F$-sequence condition from [6]. Namely,  partial orders $\leq_P $  and $\leq_Q$ must be  equivalent on corresponding  $R \subseteq P$  and $R'\subseteq Q$ isomorphic sub-posets $R$ and $R'$ which we then shall consider as identical by convention.  Note that \textbf{$R$ may not be empty!}. This requirement with respect to $P \os Q$ natural join operation is in opposition to ordinal sum  $\oplus $ which was already used in the case of antichains to produce cobweb posets [2].  For earlier,  (bipartite posets and matrices) special cases definitions  of natural join  - see [2,3,4,5,6].\\
The present overall definition of natural join $P \os Q$ of posets offers many conveniences - for example one arrives at a very simple proof of the M{\"{o}}bius function formula for cobweb posets.\\
The main aims of this article are to present at first the authors update applications of natural join and then one summarizes more systematically the general properties of natural join of posets with posing some questions arising on the way.\\ 
Let us now come over to the detailed presentation. If the posets $\left\langle P,\leq_P \right\rangle$ and $\left\langle Q,\leq_Q \right\rangle$ are disjoint then cardinal sum of posets  seams to mimic  their natural join - but these operations also in this case are different  see Fig.13. which  illustrates our choice.

\vspace{0.2cm}

\noindent \textbf{1.2.5.} Let us consider some special cases before the general definition.

\vspace{0.1cm} 

\noindent Let  posets $\left\langle P,\leq_P \right\rangle$ and $\left\langle Q,\leq_Q \right\rangle$ be \textbf{not disjoint} or constitute posets with partial orders
$\leq_P $  and $\leq_Q$ being equivalent on corresponding in $P$ and $Q$ isomorphic sub-posets  which we shall consider identical by convention. We shall say that these two posets satisfy the \textcolor{blue}{\textbf{Natural Join Conformity Condition}}.  Then we define what follows.

\vspace{0.1cm}

\noindent \textbf{1.2.5.a.} Let $P =P_1 \oplus P_2$  and  $Q = P_2 \oplus P_3$ then we define $\os$ via identity
\vspace{0.1cm}
\begin{center}
$P \os Q  \equiv \Pi_1 \oplus \Pi_2 \oplus \Pi_3 $
\end{center}

\vspace{0.1cm}
\noindent \textbf{1.2.5.b.} Let $P =P_1 + P_2$  and  $Q = P_2 + P_3$ then we define $\os$ via identity

\begin{center}
$P \os Q  \equiv P_1 + P_2 + P_3 $
\end{center}

\vspace{0.1cm}
\noindent \textbf{1.2.5.c.} Let $P =P_1 + P_2$  and  $Q = P_2 \oplus P_3$ then we define $\os$ via identity
\vspace{0.1cm}
\begin{center}
$P \os Q  \equiv P_1 + P_2 \oplus P_3 $
\end{center}

\vspace{0.2cm}

\noindent \textbf{1.2.5.d.} 

\begin {defn} {Definition}
Let posets  $P$  and $Q$  satisfy  the \textcolor{blue}{\textbf{Natural Join Conformity Condition}}. Namely - let   $\left\langle P_1, \leq_1 \right\rangle$ , $\left\langle P_2, \leq_2\right\rangle$,  $\left\langle P_3, \leq_3\right\rangle$   posets be given such that  $P =P_1 \cup P_2$  and  $Q = P_2\cup P_3$ are disjoint sums i.e. these sums represent corresponding two block partitions of $P$  and  $Q$.  Let $x \leq_1y \equiv  x \leq_2y$  for  $x,y \in P_2$.  Then we define $\os$ via identity

\begin{center}
$P \os Q  \equiv P_1 \cup P_2 \cup P_3 $
\end{center}

\vspace{0.1cm}
\noindent where in  resulting poset   $\left\langle P \os Q , \leq_{P \os Q } \right\rangle$ we put for $x,y \in P \os Q $ :

$$x \leq_{P \os Q }y  = x \leq_1y, \ x,y \in P,$$

$$x \leq_{P \os Q }y  = x \leq_2y, \ x,y \in Q.$$
The case  $x \in P_1$  and $y \in P_3$ is operated  due to transitivity of partial orders, therefore in this case

$$x \leq_{P \os Q }y  \Leftrightarrow \; \exists \; directed \; path \ x\prec\cdot p_1 \prec\cdot p_2 ...\prec\cdot p_k \prec\cdot y $$
where cover relation $\prec\cdot$  denotes either cover relation  $\prec\cdot_P$ in $P$  or cover relation  $\prec\cdot_Q$ in $Q$ depending on 
the stage of journey - if at all - along eventual directed path in Hasse DAG - DIAGram of the poset $\left\langle P \os Q,\leq_{P \os Q }\right\rangle$.
\end{defn}
  
\vspace{0.1cm}

\noindent This is to be continued in Section 4. where we shall discuss several elementary properties of natural join and other  operations on posets.  See also Figures 13,14,15.

\vspace{0.2cm}

\noindent \textbf{1.2.6.} Cobweb  posets were  identified till now with  the natural join of chain of bipartite complete digraphs  ["di-bicliques"] $B_k = \Phi_k \oplus \Phi_{k+1}$ which are ordinal sums of antichains  (see Fig.1.,Fig.3, Fig.4, Fig.5, Fig.7 and contact [2-6] and references therein) and may be therefore defined equivalently as the  ordinal sum  [linear sum ] of chain of trivially ordered sets i.e. antichains  - see Definition 4.  where cobweb  poset  $\Pi$  appears as a  linear sum of trivially ordered sets   $ \left\{\Phi_k \right\}_{k\geq0}$.

\vspace{0.1cm}

\begin{center}
 $\Pi = \oplus_{k\geq 0} \Phi_k$
 \end{center} \vspace{0.1cm}
\noindent $\left\{\Phi_k \right\}_{k\geq0}$ are then antichains i.e.  independent sets of $\Pi$. See Fig.4. and Fig.5.

\vspace{0.1cm}

\noindent Nevertheless,  Definitions 5,6,7 of  cobweb poset as \textit{natural join of relations} {binary, ternary, etc.) as in [6,5,4,3,2]  provide additionally  not only visual advantages based on sight interpretation. Namely: let $P_n$ be any $(n+1)$ - level graded poset where $n \in N \cup \left\{\infty \right\}$. Let $min(P_n)$ denotes the $P_n$ posets' set of minimal elements which we assume finite.  Let  \textit{internal} vertex $v$ means that $v\notin min(P_n)$ and the height $h(v)$ is less then the height $h(P_n)$,  $h(v)< h(P_n)=n$; $n \in N \cup \left\{\infty \right\}$.  Let us call the internal  vertex $m \in P_n$  \textbf{a \textcolor{blue}{mute} vertex} if  vertex $m$ has  either in-degree or out-degree equal zero. Let us call dummy the node  $d$ for which $in-deg(d) = 0 = out-deg(d)$, (consult Fig.2 , Fig.17 at the end of the paper and compare with Fig.1). With these naming  established we may express our crucial observation as follows.

\begin{center}
\textit{All graded posets with no mute and no dummy vertices} in their Hasse diagrams  are natural join of chain of relations and may be at the same time  interpreted an $n-ary$, $n \in N \cup \left\{\infty \right\}$ relation. 
\noindent Equivalently - \textit{zero columns or rows in bi-adjacency matrices of bipartite natural join summands of $P_n$ are forbidden}.  See figures 2 and 17 at the end of the paper and compare with other relevant.([6], [2]). 
\noindent In the case the $F(n)$-denominated $P_n$-poset,  $n \in N \cup \left\{\infty \right\}$ is infinite then $P_\infty = P$ is a natural join of chains of relations of the might be varying or constant $r$-arity ; $r \geq 2$.   
 
\end{center}  

\vspace{0.1cm}
\noindent As for cobweb posets - these have also discrete hyper-boxes representation [12,2,3, 49,17]. 

\vspace{0.1cm}

%%%%%%%%%%%%%%%%%%%%%%%%%%%%%%%%%%%%%%%%%%%%%%%%%%%%%%%%%%%%%%%%%%%%%%%%%%%%%%%%%%%%%%
\begin{figure}[ht]
\begin{center}
	\includegraphics[width=100mm]{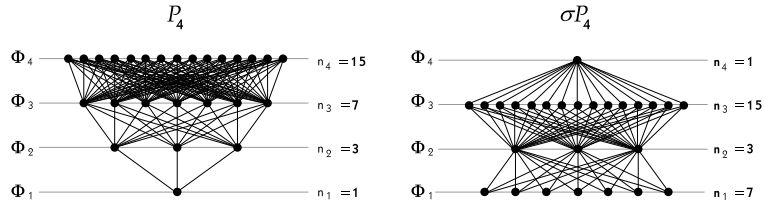}
	\caption {Display of the layer  $\langle\Phi_1 \to \Phi_4 \rangle$ = the subposet $\Pi_4$ of the  $F$ = Gaussian integers sequence $(q=2)$ $F$-cobweb poset and $\sigma \Pi_4$ subposet of the $\sigma$ permuted Gaussian $(q=2)$ $F$-cobweb poset.}\label{fig:representation}
\end{center}
\end{figure}

\begin{figure}[ht]
\begin{center}
	\includegraphics[width=100mm]{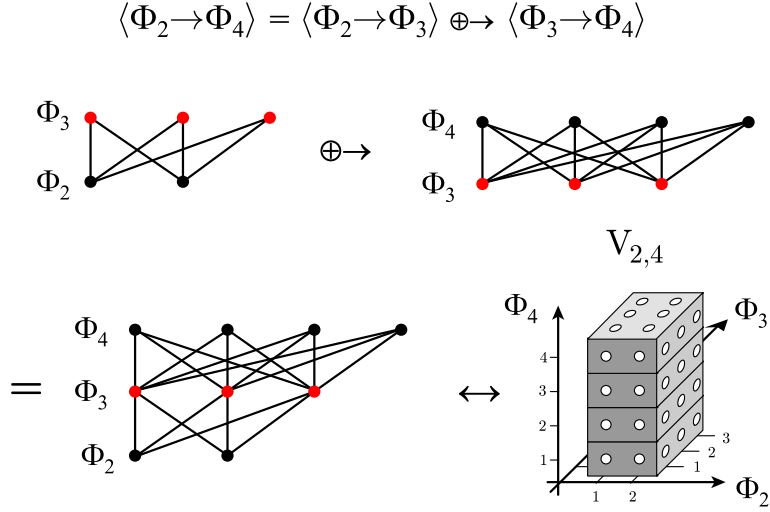}
	\caption 	{Display of the natural join of bipartite layers $\left\langle \Phi_k \rightarrow \Phi_{k+1} \right\rangle$ $F = N $, resulting in  $2\cdot3\cdot4$ maximal chains and equivalent hyper-box $V_{2,4}$ with $2\cdot3\cdot4$ white circle-dots.}\label{fig:representation} 
\end{center}
\end{figure}

\begin{figure}[ht]
\begin{center}
	\includegraphics[width=100mm]{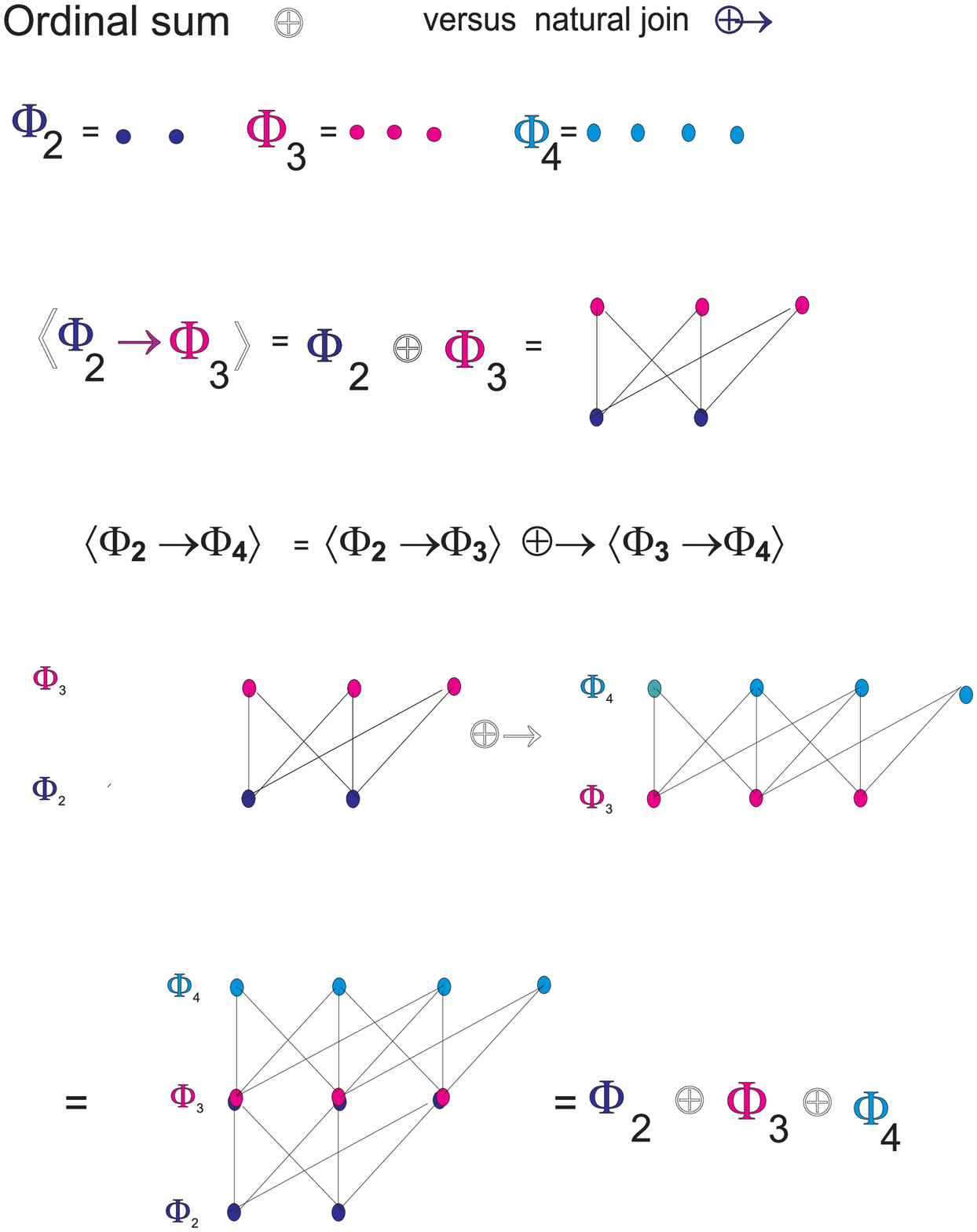}
	\caption {Display of the ordinal sum versus natural join for $F=N$.}\label{fig:representation} 
\end{center}
\end{figure}
\begin{figure}[ht]
\begin{center}
	\includegraphics[width=70mm]{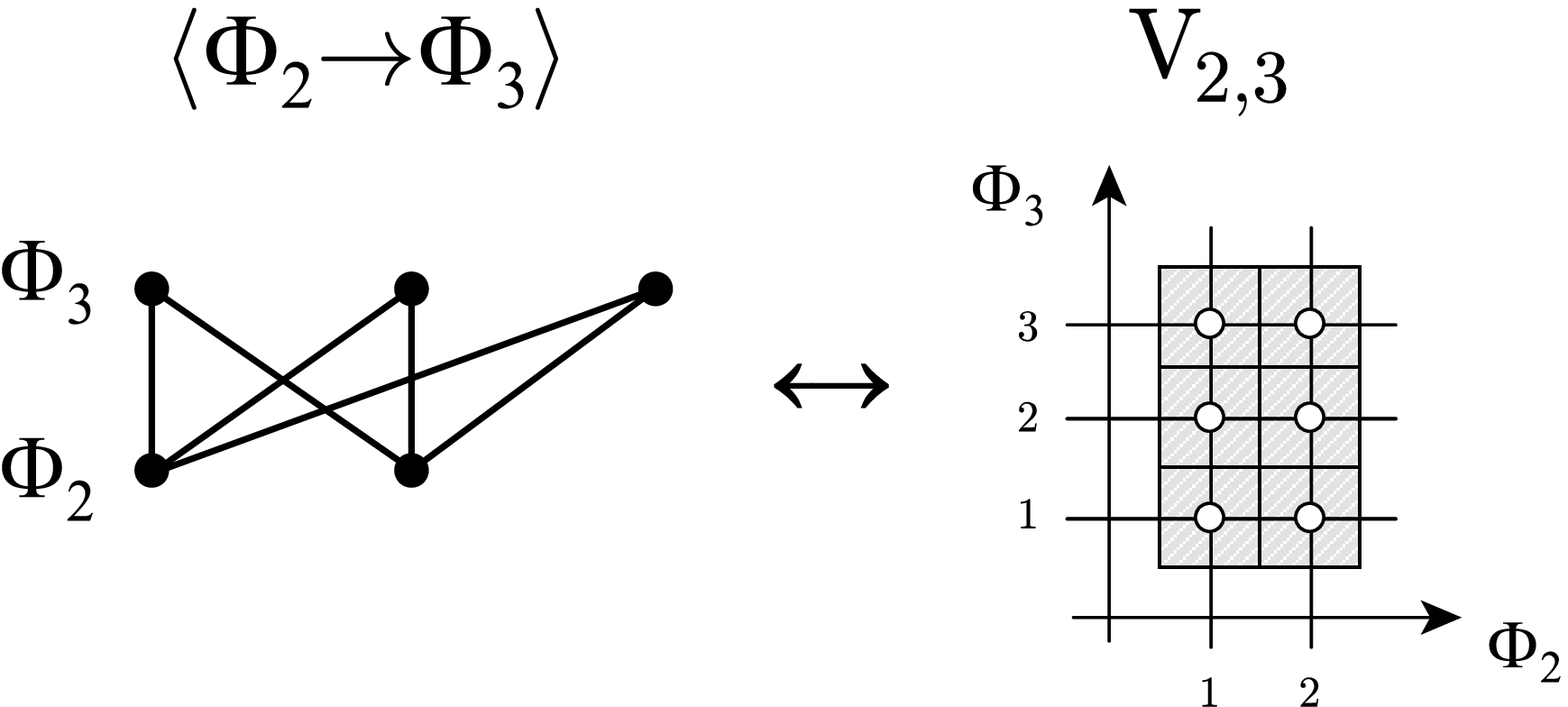}
	\caption{Bipartite  layer $\langle\Phi_3 \rightarrow \Phi_4 \rangle$ with six maximal chains and equivalent hyper-box $V_{2,3}$  with six white circle-dots} \label{fig:representation}
\end{center}
\end{figure}
%%%%%%%%%%%%%%%%%%%%%%%%%%%%%%%%%%%%%%%%%%%%%%%%%%%%%%%%%%%%%%%%%%

\vspace{0.1cm}
\noindent Illustration  -  see Fig.4, Fig.8, Fig.6, Fig.7.

\begin{figure}[ht]
\begin{center}
	\includegraphics[width=70mm]{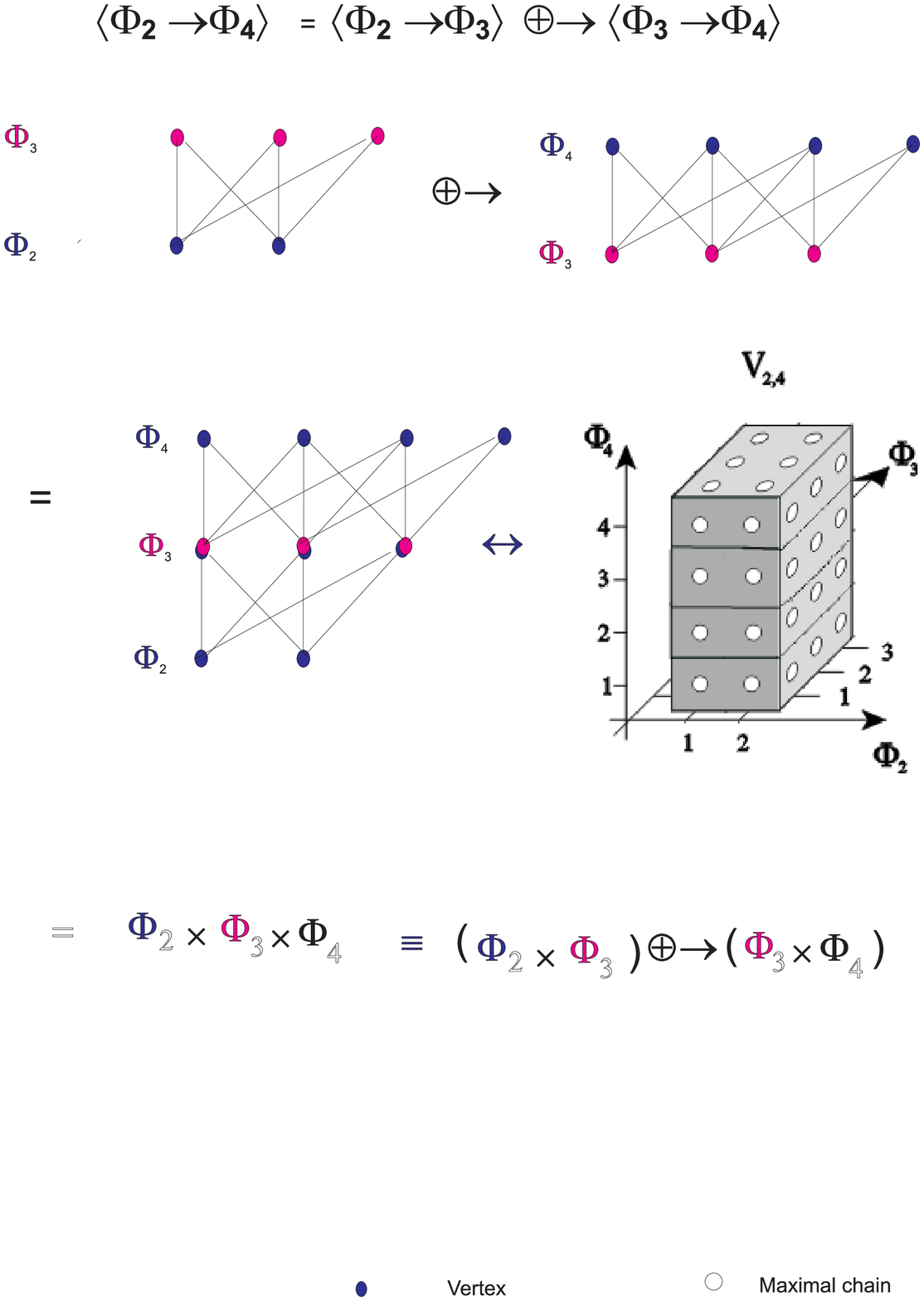}
	\caption{Natural join of layers versus Cartesian product constituting the hyper-box $V_{2,4}$.} \label{fig:representation}
\end{center}
\end{figure}
\begin{figure}[ht]
\begin{center}
	\includegraphics[width=70mm]{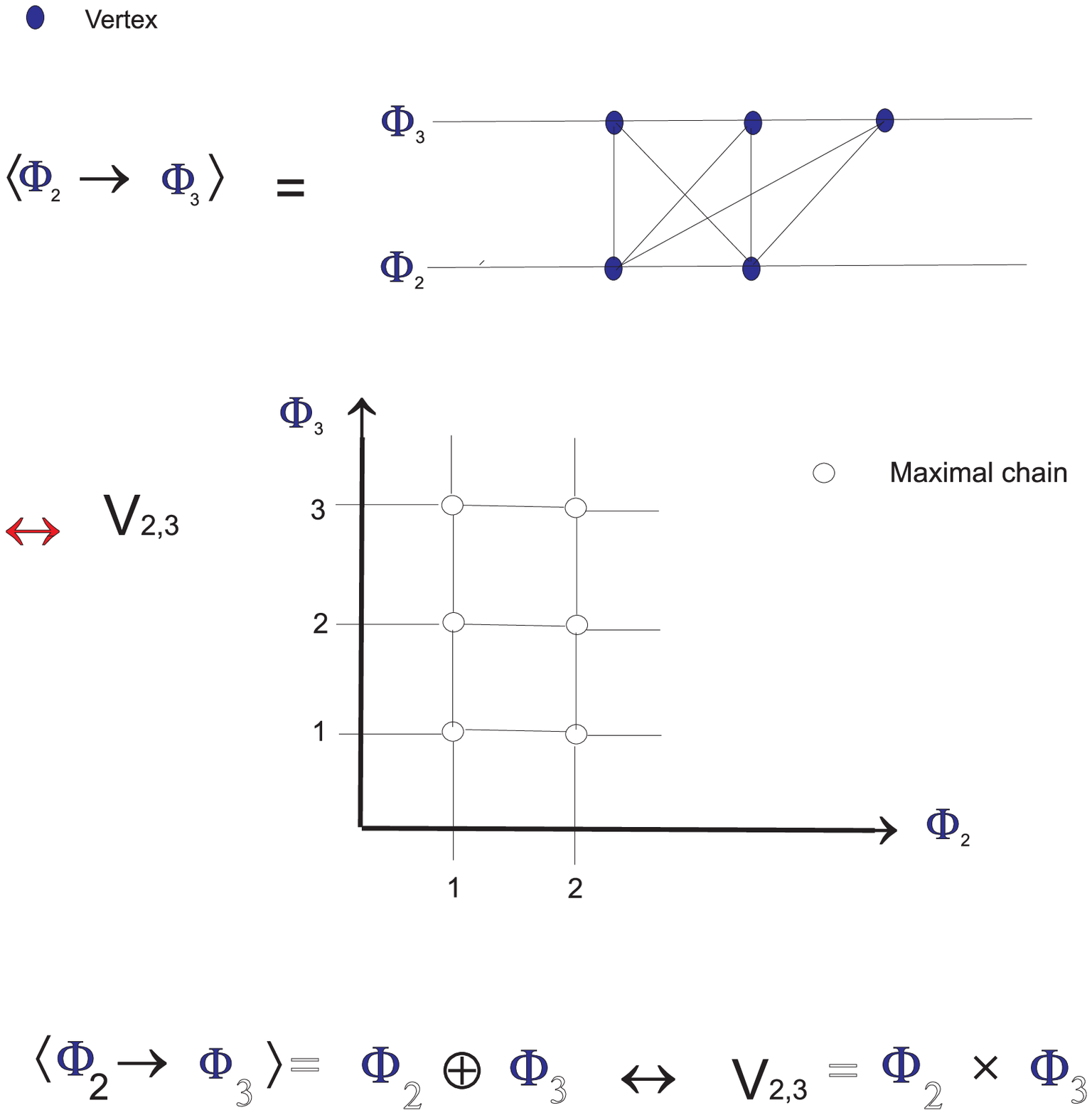}
	\caption{Natural join as ordinal sum of levels $\Phi_2\os \Phi_3 \equiv  \Phi_2\oplus \Phi_3 \equiv V_{2,3}$.}\label{fig:representation}
\end{center}
\end{figure}

%%%%%%%%%%%%%%%%%%%%%%%%%%%%%%%%%%%%%%%%%%%%%%%%%%%%%%%%%%%%%%%%%%%%%%%%%%%%%%%%%%%%%%

\vspace{0.1cm}

\noindent \textcolor{blue}{\textbf{Heralding influence posets.}} Apart  from definitions below in terms of natural join of $r$-ary relations' sequence  (i.e. we deal with ordered partition of cobweb poset  into  $r$-level layers) - apart from these we may interpret  cobweb poset  as the special extreme influence poset [53,54]. In more detail.

\vspace{0.1cm}
%%%%%%%%%%%%%%%%%%%%%%%%%%%%%%%%%%%%%%%%%
\noindent \textbf{1-ary.} Consider a chain  $\left\{\Phi_k\right\}_{k\geq 0}$  of  trivial posets  $\equiv$ antichains  $\equiv$ a chain of  independent sets  $\equiv$ a  chain of  trivial unary relations, then

\begin{defn} (cobweb poset $\Pi_n$) 
\noindent Let  $n\in N \cup \left\{0\right\}\cup \left\{\infty\right\}$. Then

$$\Pi_n =  \oplus_{s= 0}^n \Phi_s ,$$

$$\Pi =  \oplus_{s\geq 0} \Phi_s, $$
where $\oplus$ denotes ordinal sum of posets  hence

$$\langle\Phi_k \to \Phi_n \rangle \equiv \oplus_{n\geq s\geq k} \Phi_s $$
is called the layer of the cobweb poset  [2-12] while $\left\langle \Phi_k\right\rangle_{k\geq 0}$    
constitute the sequence  of independent sets  i.e.  the sequence  of antichains . 
\end{defn}

\noindent\textbf{2-ary.} Consider the chain  $\left\{B_k\right\}_{k\geq 0}$  of di-bicliques   $B_k= \Phi_k \oplus \Phi_{k+1} $  $\equiv$ the chain of universal binary relations. then

\begin{defn} (cobweb poset $\Pi$)

$$\Pi = \os_{k \geq 0}B_k $$
where $\os $ denotes natural join  of posets.  Here   $\left\langle B_k\right\rangle_{k\geq 0}$    
constitute the sequence  of  universal binary relations $\equiv$  the sequence of di-bicliques.   
\end{defn}

\noindent \textbf{3-ary.} Consider the chain  $\left\{T_k\right\}_{k\geq 0}$  of  $T_k= \Phi_k \oplus \Phi_{k+1} \oplus \Phi_{k+2}$    i.e  the   chain of ternary universal relations, then

\begin{defn} (cobweb poset $\Pi$)
$$\Pi =  \os _{k\geq 0}T_k $$
where $\os $ denotes natural join  of posets.  Here   $\left\langle T_k\right\rangle_{k\geq 0}$    
constitute the sequence  of  ternary universal relations.
\end{defn}

\noindent \textbf{r-ary.} Finally, consider the chain  $\left\{R_k\right\}_{k\geq 0}$  of  $R_k= \Phi_k \oplus \Phi_{k+1}\oplus... \oplus \Phi_{k+r-1}$ , $r>0$    i.e  the  chain of $r$-nary universal relations, then

\begin{defn} (cobweb poset $\Pi$)
$$\Pi =  \os _{k\geq 0} R_k $$
where $\os $ denotes natural join  of posets .  Here   $\left\langle R_k\right\rangle_{k\geq 0}$    
constitute the sequence  of  $r$-nary universal relations.
\end{defn}

\noindent \textbf{Comment 1.} Bring together and make \textbf{identifications} of graded KoDAGs\textbf{ with  $n$-ary relations} digraph representation as in [2,3,4,5,6]:

$$ \leq = \Phi_0\times\Phi_1\times ... \times\Phi_n \Longleftrightarrow cobweb \  poset  \Longleftrightarrow  KoDAG ,$$
for the natural join of di-bicliques and similarly for $\leq$ being  natural join  of any sequence binary relations [see Fig.1,3,4,5,7,11]

$$ \leq \subseteq \Phi_0\times\Phi_1\times ... \times\Phi_n \Longleftrightarrow cobweb \  poset  \Longleftrightarrow  KoDAG .$$
\textcolor{red}{\textbf{Warning.}} Note that \textcolor{red}{\textbf{not for all}} $F$-graded posets their partial orders may be consequently identified with 
$n$-ary relations, where  $F = \left\langle k_F \right\rangle_{k=1}^n $ while $n \in N \cup \left\{\infty\right\}$. This is possible iff
no  biadjacency matrices entering the natural join for $\leq$  has a zero column or a zero row. Recall if  an internal  vertex $m \in P$ has not 
either incoming or outgoing arcs then we shall call it the \textcolor{red}{\textbf{mute}} node [2]. Recall.Dummy vertex is that  $d \in P$ for which $in-deg(d) = 0 = in-deg(d)$.  This naming being adopted we may say now:

\vspace{0.1cm}

\noindent $F$-\textit{graded poset may be identified with} $n$-\textit{ary relation  as above iff it is} $F$-\textit{graded poset with neither mute no dummy nodes.}\\ 
\vspace{0.1cm} 
\noindent Equivalently - \textit{zero columns or rows in bi-adjacency matrices of bipartite natural join summands of $P$ are forbidden}.  See and compare with relevant figures 2 and 17 at the end of the paper and perhaps consult [6,2].

\vspace{0.1cm}
\noindent  \textcolor{blue}{\textbf{Heralding.}} Apart  from definitions above in terms of natural join of $r$-ary relations' sequence  (i.e. we deal with ordered partition of cobweb poset  into  $r$-level layers) - apart from these we may interpret  cobweb poset  as the special extreme influence poset [53,54].
%%%%%%%%%%%%%%%%%%%%%%%%%%%%%%%%%%%%%%%%%%%%%%
%%%%%%%%%%%%%%%%%%%%%%%%%%%%%%%%%%%%%%%%%%%%%%

\noindent This is to be reported on next in the  subsection  2.2.  Meanwhile let us continue with what follows.\\

\vspace{0.2cm}

%%%%%%%%%%%%%%%%%%%%%%%%%%%%%%%%%%%%%%%%%%%%%%%%%%%%%%%%%%%%%%%%%%
%%%%%%%%%%%%%%%%%%%%%%%%%%%%%%%%%%%%%%%%%%%%%%%%%%%%%%%%%%%%%%%%%%

%%%%%%%%%%%%%%%%%%%%%%%%%%%%%%%%%%%%%%%%%%%%%%%%%%%%%%%%%%%%%%%%%%%%%%%%%%%%%%%%%%%%%%%%%%%%
\noindent \textbf{1.3. Convene  the above with hyper-boxes from [12,2]} \\

%%%%%%%%%%%%%%%%%%%%%%%%%%%%%%%%%%%%%%%%%%%%%%%%%%%%%%%%%%%%%%%%%%
%%%%%%%%%%%%%%%%%%%%%%%%%%%%%%%%%%%%%%%%%%%%%%%%%%%%%%%%%%%%%%%%%%

\noindent  Recall [3,9]: $C_{max}(\Pi_n)$ is the set of all maximal chains of $\Pi_n$. Recall [3,9]: 
$$C^{k,n}_{max} = \big\{ \mathrm{maximal\ chains\ in\ } \langle \Phi_k \rightarrow \Phi_n \rangle \big\}.$$

\noindent Consult  Section 3. in  [12] in order to view $C_{max}(\Pi_n)$ or  $C^{k,n}_{max}$ as the hyper-box of points.\\ 
Namely [12,3,2] denoting with $V_{k,n}$ the discrete finite rectangular $F$-hyper-box or $(k,n)-F$-hyper-box or in everyday parlance just $(k,n)$-box
$$
	V_{k,n} = [k_F]\times [(k+1)_F]\times ... \times[n_F]
$$
\noindent we identify (see Figure 7.) the following two just by agreement according to the $F$-natural identification:
$$
	C^{k,n}_{max} \equiv V_{k,n}
$$
i.e.
$$
C^{k,n}_{max} = \big\{ \mathrm{maximal\ chains\ in\ } \langle \Phi_k \rightarrow \Phi_n \rangle \big\} \equiv V_{k,n}.
$$
%%%%%%%%%%%%%%%%%%%%%%%%%%%%%%%%%%%%%%%%%%%%%%
\vspace{0.1cm}
\noindent Illustration  -  see Fig.4, Fig.6, Fig.7, Fig.8  and for more  see [12], [48], [12] and  [3]  with specific indication on [17].

\vspace{0.2cm}

\noindent \textcolor{red}{\textbf{Important}}.  Accordingly  the \textcolor{blue}{\textbf{natural join operation}} of  discrete hyper-boxes - ( cobweb posets are  encoded by discrete hyper-boxes [12])  is just \textcolor{blue}{\textbf{Cartesian product of them accompanied with projection out}} of sine qua non common faces (see \textcolor{blue}{\textbf{Fig.7}}) which is schematically represented by an sample case below  [note - cobweb poset might be defined also as the ordinal sum of its independence sets (antichains) $\left\{\Phi_k\right\}_{k \geq )}$ ],  
\vspace{0.1cm}
$$\left(\Phi_k \times \Phi_{k+1}\right) \os \left(\Phi_{k+1} \times \Phi_{k+2}\right) \equiv  \Phi_k \times \Phi_{k+1} \times \Phi_{k+2}, $$
see Fig.6 , Fig.7, Fig.8.

\vspace{0.2cm}

\noindent \textcolor{red}{\textbf{Important}}. \textbf{The face lattice of a polytope coding} [53]. \\  
\noindent Compare the above  discrete hyper-boxes coding representation of unbounded number of $F$-cobweb posets with one Agnarsson's model example of a classical poset $F_{\pi} =$ the face lattice of an $n$-dimensional polytope $\pi$ - viewed as stacked height-$2$ sub-posets $F_{\pi}(k, k + 1)$ , one on top of the other [53], which exactly means that these bipartite subposets $F_{\pi}(k, k + 1)$ are just naturally  joined by $\os$ binary operator.  Here in Agnarsson's paper [53] the subposet  $F_{\pi}(k, k + 1)$ is the height-$2$ sub-poset  i.e. bipartite sub-poset of $F_{\pi}$ consisting of the $k$ and 
$(k + 1)$-dimensional faces of $\pi$. Consequently $F_{\pi}$  is  thought of being formed by stacking  [ i.e.  by natural  joining !]  $F_{\pi}(k, k + 1)$  on top of   $F_{\pi}(k-1, k)$  for each   $k = 0, 1, . . . , n$. Compare with Figures 1,2,3,4,8. Admit $n \in N \cup \left\{\infty\right\}$.

\noindent In the $F_{\pi}$ case of the face lattice of an $n$-dimensional polytope $\pi$  the  Agnarsson's stacking  is exactly the natural  join $\os$ operations of layers. This natural join operation is    naturally natural since $F_{\pi}$  is the graded poset equipped  with a grading function mapping each face of the polytope $\mu$ (i.e. each element of the poset $F_{\pi}$) to its dimension  [53].

\vspace{0.3cm}

%%%%%%%%%%%%%%%%%%%%%%%%%%%%%%%%%%%%%%%%%%%%
\noindent \textbf{1.4. Cobweb posets' M\"{o}bius function}\\

\noindent  In order to find out the  M\"{o}bius function $\mu$ of the cobweb poset $\Pi$ note the following obvious statements  ($\left|\Phi_k\right|= k_F$)
  
\vspace{0.1cm}
  
\noindent Obvious statement: $\mu(x,x) = 1, \   \mu(x,y) = - 1$ for $x\prec\cdot y , \ $and$  \   \mu(x,z) = k_F - 1 $  for   $[x,z] = x\oplus \Phi_k\oplus z $.

\vspace{0.1cm}

\noindent Obvious statement: ($\left|\Phi_k\right|= k_F$)  for  $[x,z] = x\oplus \Phi_k\oplus z $ , $x\in \Phi_{k-1}, \ z \in \Phi_{k+1}$

\vspace{0.1cm}
\noindent Obvious statement: $$\mu(x,z) = [ k_F - 1].$$ 
\vspace{0.1cm}
\noindent Obvious statement ($\left|\Phi_k\right|= k_F$)  :  for  $[x,z] = x\oplus \Phi_k\oplus \Phi_{k+1}\oplus z \equiv x\oplus B_k \oplus z$ , $x\in \Phi_{k-1}, \ z \in \Phi_{k+2}$ 
\vspace{0.1cm}
$$\mu(x,z) = - [ k_F - 1][ (k+1)_F - 1].$$
\vspace{0.1cm}
\noindent Illustration  ( see Fig.9 and Fig.10).

\vspace{0.1cm}  

\begin{figure}[ht]
\begin{center}
	\includegraphics[width=70mm]{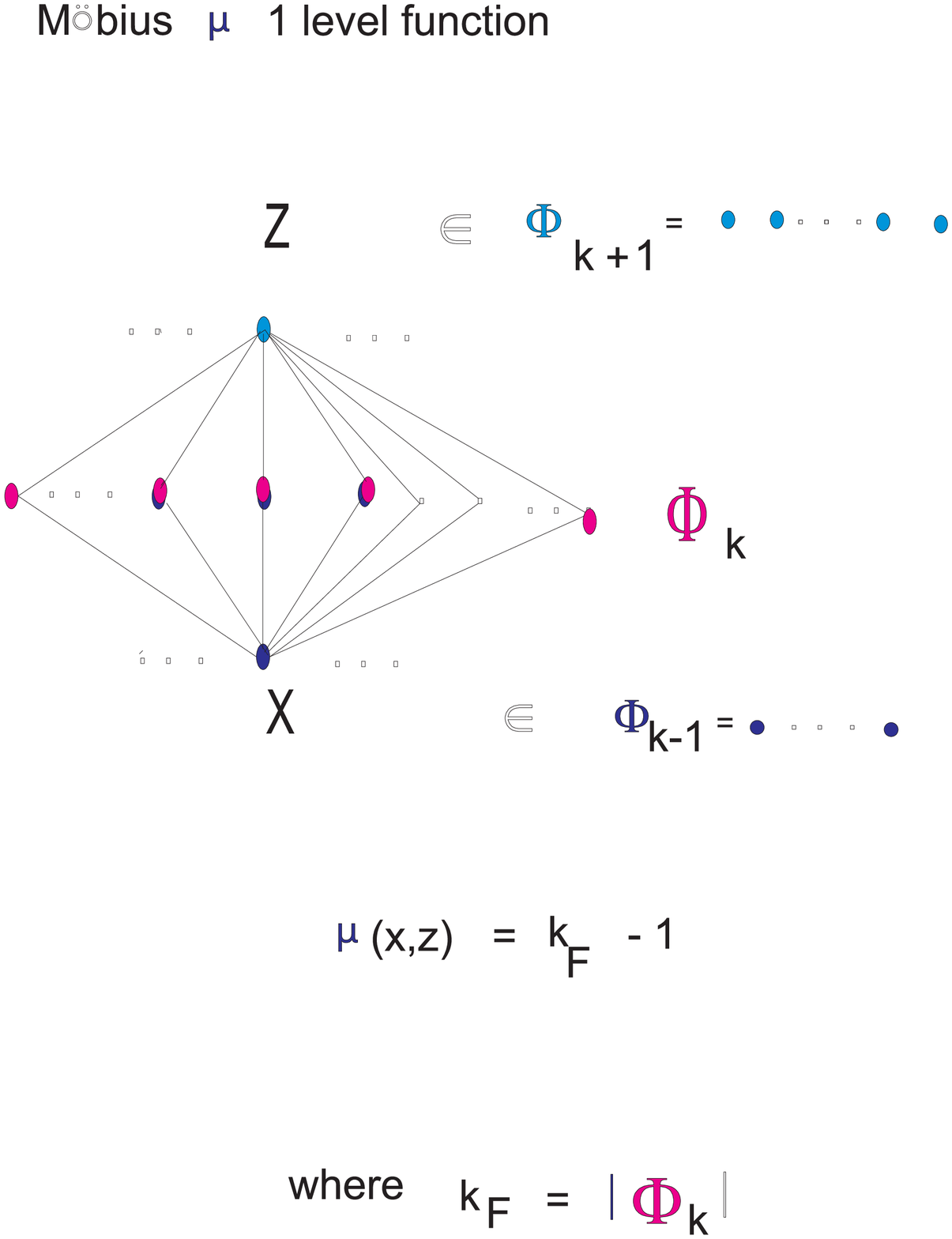}
	\caption{M\"{o}bius 1 level function for  $[x,z] = x\oplus \Phi_k\oplus z $ , $x\in \Phi_{k-1}, \ z \in \Phi_{k+1}$.}\label{fig:representation}
\end{center}
\end{figure}

\begin{figure}[ht]
\begin{center}
	\includegraphics[width=70mm]{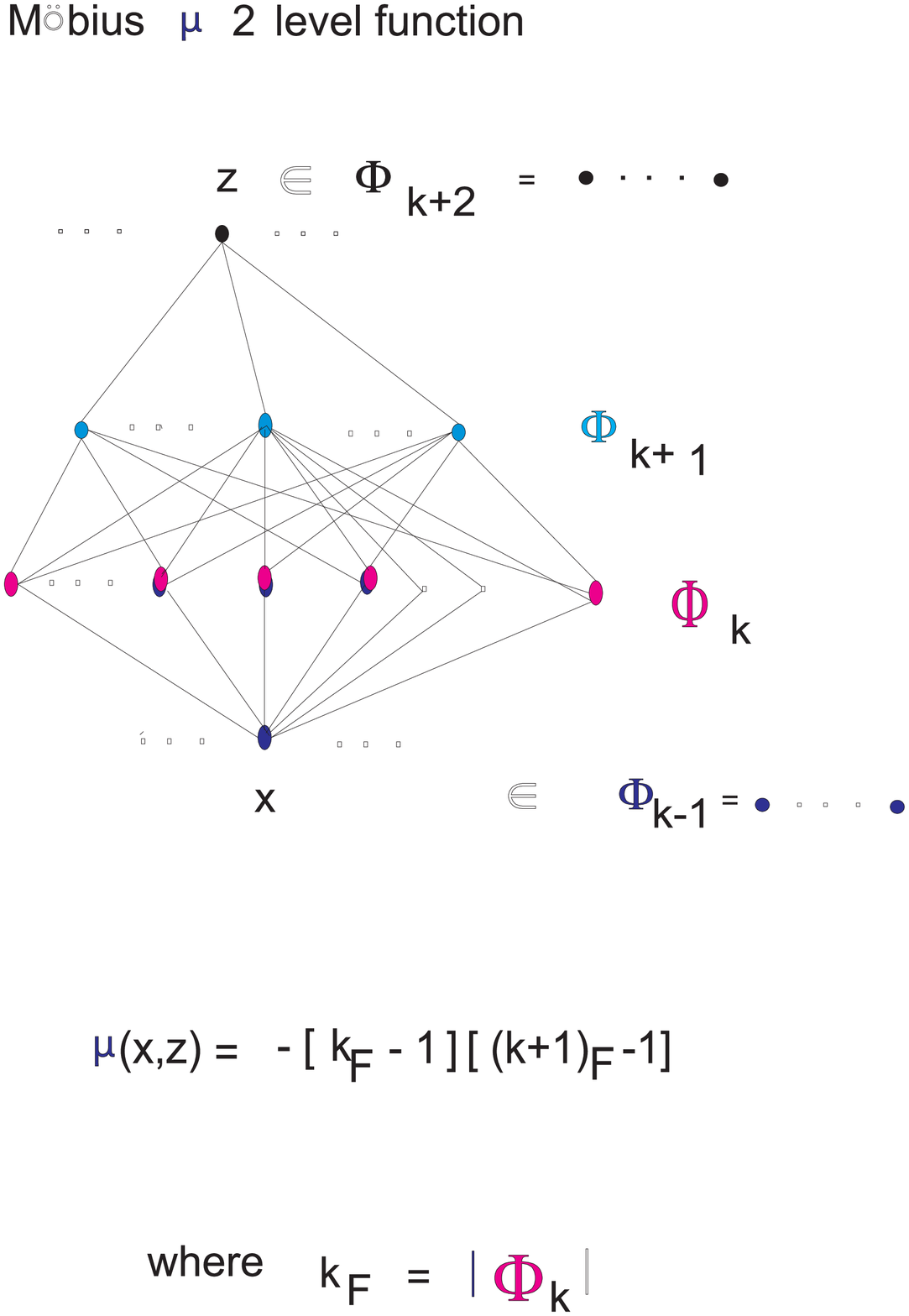}
	\caption{M\"{o}bius 2 level function  for  $[x,z] = x\oplus \Phi_k\oplus \Phi_{k+1}\oplus z $ , $x\in \Phi_{k-1}, \ z \in \Phi_{k+2}$.}\label{fig:representation}
\end{center}
\end{figure}

\vspace{0.1cm} 

\begin{figure}[ht]
\begin{center}
	\includegraphics[width=70mm]{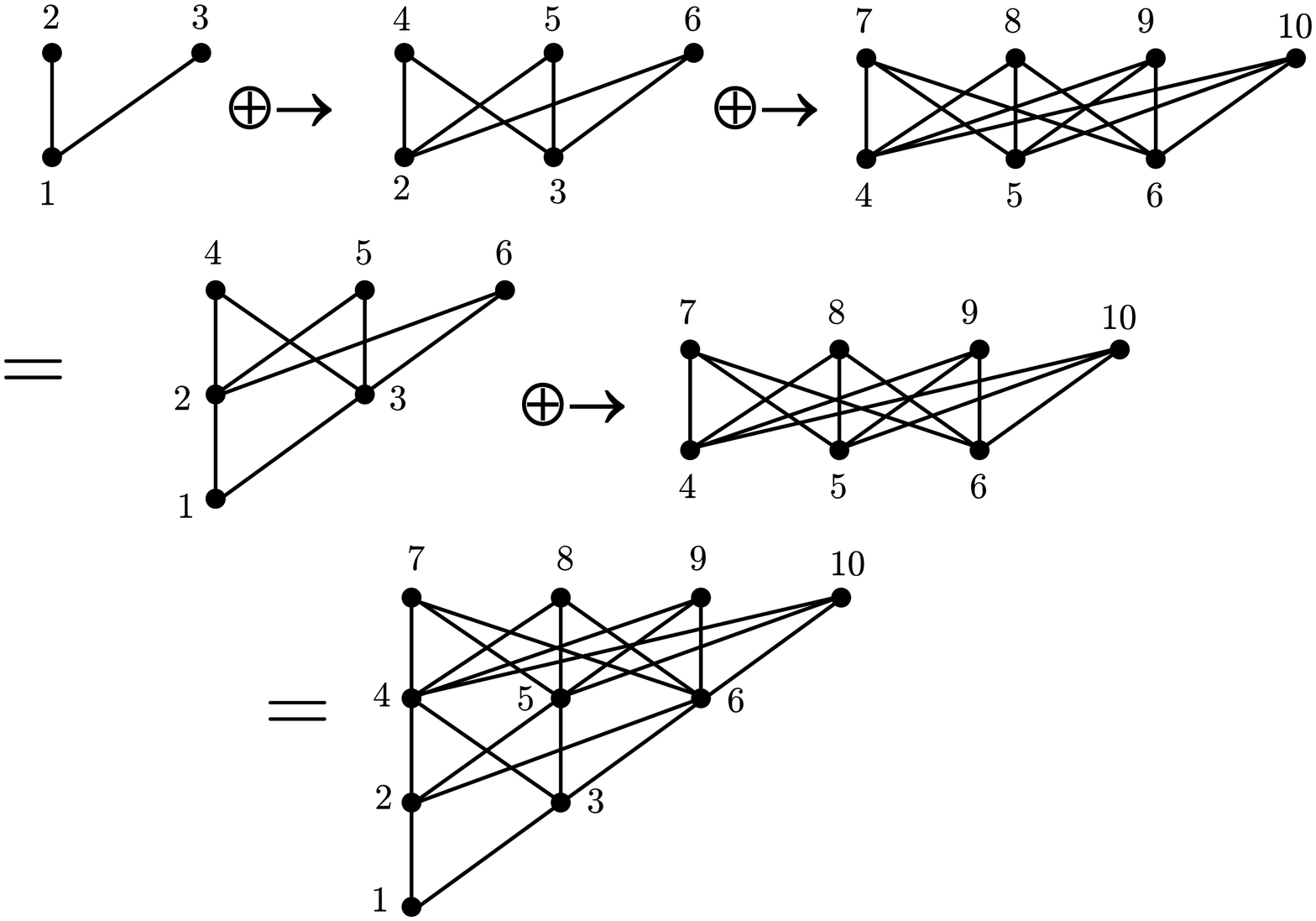}
	\caption{Natural join in natural labeling.} \label{fig:representation}
\end{center}
\end{figure}

\vspace{0.1cm} 

\noindent Hence via induction  for  $x \in \Phi_r$ ,   $z \in \Phi_s$  and  for     
$$s>r , \ [x,z] = x\oplus \Phi_{r+1}\oplus...\oplus \Phi_{s-1}\oplus z$$ 
we have the final obvious statement:
  $$\mu(x,y) = (-1)^{s-r}\prod_{k=r+1}^{s-1}[k_F-1]. $$

%%%%%%%%%%%%%%%%%%%%%%%%%%%%%%%%%%%%%%%%%%%%%%

\vspace{0.1cm} 

\noindent Compare the above  formula from [2] with the form  derived and discussed in [3] and those declared in [47,25, 29] for Fibonacci sequence. See also [27,28].
\vspace{0.2cm}
\noindent Naturally the values of $\mu(x,y)$ depend only on the rank of its arguments - here $r(x)=r$ and $r(y)=s$ - which is the reason of coding matrix existence 
for M\"{o}bius function in a natural labeling representation (see: examples in [3,2]). The rank function is here defined as follows: $r(x)= r$  if  $x\in \Phi_r$.
%%%%%%%%%%%%%%%%%%%%%%%%%%%%%%%%%%%%%%%%%%%%%%

%%%%%%%%%%%%%%%%%%%%%%%%%%%%%%%%%%%%%%%%%%%%%%

%%%%%%%%%%%%%%%%%%%%%%%%%%%%%%%%%%%%%%%%%%%%%%

\section{Combinatorial interpretation.}
\noindent \textbf{2.1. Kwasniewski combinatorial interpretation.} For \textbf{combinatorial interpretation of cobweb posets} via their cover relation digraphs (Hasse diagrams) called KoDAGs see [8,7,6,2,3,4]. The recent equivalent formulation of this combinatorial interpretation is to be found in [7] (Feb 2009) or [9] from which we quote it here down.

\begin{defn}
$F$-\textbf{nomial} \textbf{coefficients} are defined as follows
$$
	\fnomial{n}{k} = \frac{n_F!}{k_F!(n-k)_F!} 
	= \frac{n_F\cdot(n-1)_F\cdot ...\cdot(n-k+1)_F}{1_F\cdot 2_F\cdot ... \cdot k_F}
	= \frac{n^{\underline{k}}_F}{k_F!}
$$
\noindent while $n,k\in \mathbb{N}$ and $0_F! = n^{\underline{0}}_F = 1$  with $n^{\underline{k}}_F \equiv \frac{n_F!}{k_F!}$ staying for falling factorial.
$F$ is called  $F$-graded poset  \textcolor{blue}{\textbf{admissible}} sequence iff  $\fnomial{n}{k} \in N \cup\left\{0\right\}$ ( In particular we shall use the expression - $F$-cobweb admissible sequence).
\end{defn}

\begin{defn}
$$
C_{max}(\Pi_n) \equiv  \left\{c=<x_0,x_1,...,x_n>, \: x_s \in \Phi_s, \:s=0,...,n \right\} 
$$  
i.e. $C_{max}(\Pi_n)$ is the set of all maximal chains of $\Pi_n$
\end{defn}

\noindent Consequently (see Section 2 in [12]  on Cobweb posets' coding via $N^\infty$ lattice boxes) we introduce natural notation as follows.

\begin{defn} ($C^{k,n}_{max}$) Let  

$$
C_{max}\langle\Phi_k \to \Phi_n \rangle \equiv \left\{c=<x_k,x_{k+1},...,x_n>, \: x_s \in \Phi_s, \:s=k,...,n \right\}\equiv
$$

$$
	\equiv \big\{ \mathrm{maximal\ chains\ in\ } \langle \Phi_k \rightarrow \Phi_n \rangle \big\} \equiv
	C_{max}\big( \langle \Phi_k \rightarrow \Phi_n \rangle \big) \equiv
	C^{k,n}_{max}.
$$
\end{defn}

\noindent \textbf{Note.} The $C_{max}\langle\Phi_k \to \Phi_n \rangle \equiv C^{k,n}_{max}$
is the hyper-box points'  set [12,2] of  Hasse sub-diagram corresponding maximal chains and it defines biunivoquely 
the layer $\langle\Phi_k \to \Phi_n \rangle = \bigcup_{s=k}^n\Phi_s$  as the set of maximal chains' nodes (and vice versa) -
for  these arbitrary $F$-denominated \textbf{graded} DAGs (KoDAGs included).

\vspace{0.1cm}

\noindent The equivalent to that of [8,7,9] formulation of the fractals reminiscent combinatorial interpretation of cobweb posets via their cover relation digraphs (Hasse diagrams) is the following.

\vspace{0.2cm}

\noindent \textbf{Theorem 1} [11,9,8,3] \\
\noindent(Kwa\'sniewski) \textit{For $F$-cobweb admissible sequences $F$-nomial coefficient $\fnomial{n}{k}$ is the cardinality of the family of \emph{equipotent} to  $C_{max}(P_m)$ mutually disjoint maximal chains sets, all together \textbf{partitioning } the set of maximal chains  $C_{max}\langle\Phi_{k+1} \to \Phi_n \rangle$  of the layer   $\langle\Phi_{k+1} \to \Phi_n \rangle$, where $m=n-k$.}

\vspace{0.1cm} 
\noindent For environment needed and then  simple combinatorial proof see [8,9,4,5]  easily accessible via Arxiv.

\vspace{0.1cm}

\noindent \textbf{Comment 2}. For the  above Kwa\'sniewski combinatorial  interpretation of  $F$-nomials' array \textit{it does not matter}  of course whether the diagram is being directed  or not, as this combinatorial interpretation is  equally valid for partitions  of the family of  $SimplePath_{max}(\Phi_k - \Phi_n)$ in  comparability graph of the Hasse  digraph with self-explanatory notation used on the way. The other insight into this irrelevance for combinatoric interpretation is [9]: colligate the coding of $C^{k,n}_{max}$ by hyper-boxes. (More on that soon).  And to this end recall  what really also matters here : a poset is graded if and only if every connected component of its \textbf{comparability graph} is graded. We are concerned here with connected graded graphs and digraphs.

\vspace{0.1cm}

\noindent For the relevant recent developments see [10]  while \textbf{[11]} is their all \textbf{source paper} as well as those reporting on the broader affiliated research (see [12-23,25-29,47-49] and references therein). The inspiration for  "`philosophy"'  of notation in mathematics as that in Knuth's from [24] - in the case of "`upside-downs"'  has been  driven by Gauss "`$q$-Natural numbers"'$\equiv N_q = \left\{n_q=q^0+q^1+...+q^{n-1}\right\}_{n\geq 0}$ from finite geometries of  linear subspaces lattices over Galois fields. As for the earlier use and origins of the use of this author's upside down notation see [30-46].

\vspace{0.1cm}

\noindent In discrete hyper-boxes language the combinatorial interpretation reads: 

\vspace{0.2cm}

\noindent \textbf{Theorem 2}\\
\noindent \textit{For $F$-cobweb admissible sequences $F$-nomial coefficient} $\fnomial{n}{k}$ \textit{is the cardinality of the family of \emph{equipotent} to}  $V_{0,m}$ \textit{mutually disjoint 
discrete hyper-boxes, all together \textbf{partitioning } the discrete hyper-box }  $V_{k+1,n}$    $\equiv$   \textit{the layer}   $\langle\Phi_{k+1} \to \Phi_n \rangle$,\textit{ where} $m=n-k$.

\vspace{0.2cm}

\noindent \textbf{Comment 3.} General "environmental"  comment. Recall:  all graded posets with no mute vertices in their Hasse diagrams [3,2] (i.e. no vertex has  in-degree or out-degree equal zero) are natural join of chain of relations and may be at the same time  interpreted an $n-ary$, $n \in N \cup \left\{\infty \right\}$ relation.\\
\noindent Scrape together  any binary relation $R$ with Hasse digraph cover relation  $\prec\cdot$  and identify as in [4,3] $\zeta({\bf R})\equiv {\bf
R}^{*}$ with incidence algebra zeta function and with zeta matrix of the poset associated to its Hasse digraph, where the \textbf{reflexive} reachability  relation $\zeta({\bf R})\equiv {\bf
R}^{*}$ is defined as

$${\bf R^*}=R^0\cup R^1\cup R^2\cup\ldots\cup R^n\cup\ldots  
\bigcup_{k>0}{R}^k={\bf R}^{\infty}\cup {\bf I}_A =$$
 
\begin{center}
=  {\bf transitive} and {\bf reflexive} closure of $\bf{R}$  $\Leftrightarrow$
\end{center}
 
$$
\Leftrightarrow\;  A(R^{\infty})=A({ R})^{\copyright 0}\vee A({
R})^{\copyright 1}\vee A({ R})^{\copyright 2}\vee \ldots \vee A(R)^{\copyright n}\vee \ldots ,
$$
where $A({\bf R})$ is the Boolean  adjacency matrix of the
relation ${\bf R}$ simple digraph  and  $\copyright$ stays for Boolean product.

\vspace{0.1cm}

\noindent Now collocate  and/or recall from [3] the resulting schemes.
\noindent {\bf Schemes:}
 $$ <=\prec\cdot ^{\infty}=\;connectivity\;of\;\prec\cdot$$
 $$ \leq=\prec\cdot ^{*}=\;reflexive \; reachability\;of\;\prec\cdot$$
 $$\prec\cdot ^{*}=\zeta(\prec\cdot).$$ 
 
\vspace{0.1cm}

%%%%%%%%%%%%%%%%%%%%%%%%%%%%%%%%%%%%%%%%%%%%%%%%%%%%%%%%%%%%%%%%%%%%%%%%%%%%%%%%%
\noindent \textbf{Comment 4.} General  "fractal-reminiscent" comment. The discrete $m$-dimensional $F$-box ($m = n-k$) with edges' sizes designated by natural numbers' valued sequence $F$ 
where invented in [12] as a response to the so called \emph{Kwa\'sniewski cobweb tiling} problem posed in [11] and then repeated in [8]. This tiling problem was considered by Dziemia\'nczuk 
in [13] where it was shown that not all admissible $F$-sequences permit tiling as defined in [11]. Then this tiling  problem (see Fig.12 from [48]) was reformulated by Dziemia\'nczuk in discrete hyper-boxes language [48] in two steps as follows. The first step: the now simple observation concerning Kwa\'sniewski  combinatorial interpretation was expressed in hyper-boxes language.

\vspace{0.1cm}

\begin{figure}[ht]
\begin{center}
	\includegraphics[width=100mm]{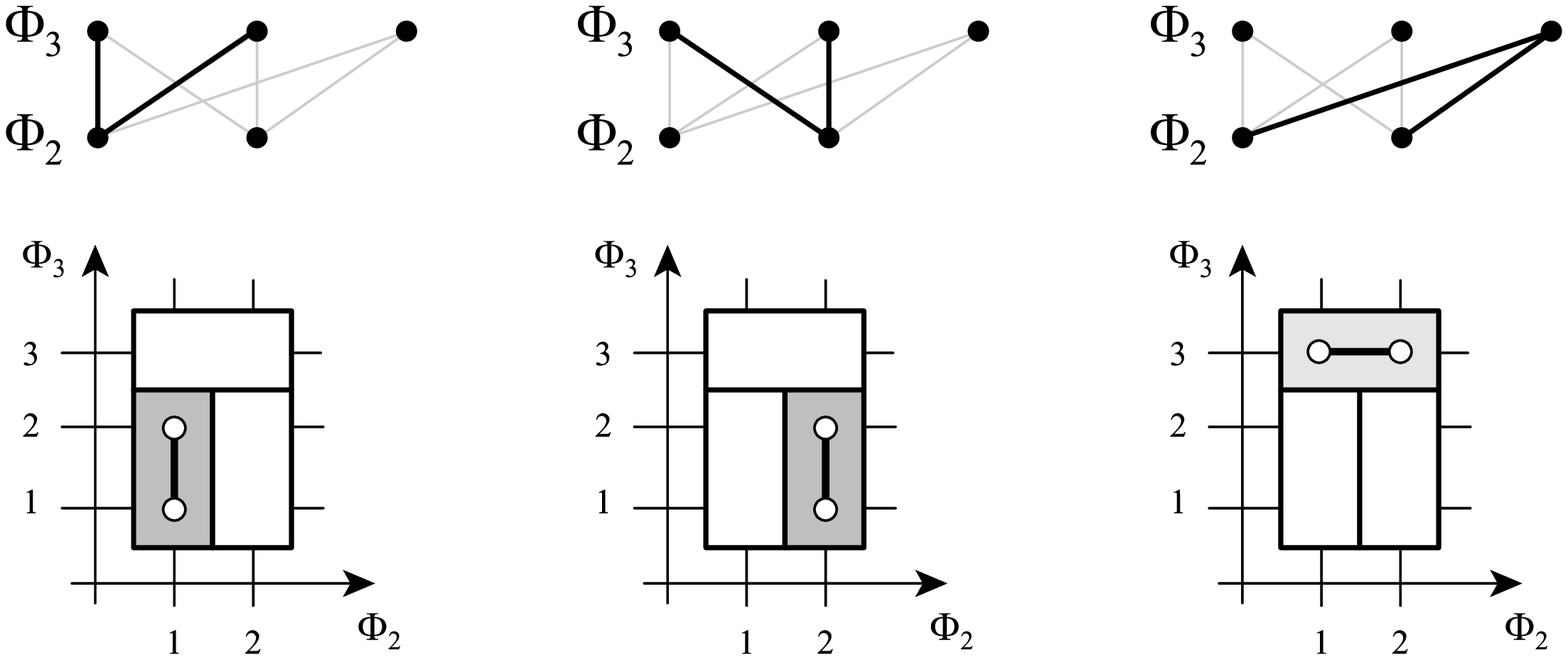}
	\caption{Correspondence between tiling of $F$-box $V_{3,4}$ with two white-dots boxes and tiling of the $\layer{3}{4}$ with two-chain subbposets.}
\end{center}
\end{figure}

\noindent \textbf{Fact} (Kwa\'sniewski [11,8])\\
Let $F$ be an admissible sequence. Take any natural numbers $n,m$ such that $n\geq m$, then the value of $F$-nomial coefficient $\fnomial{n}{k}$ is equal to the number of sub-boxes that constitute a  $\kappa$-partition of $m$-dimensional $F$-box $V_{m,n}$ where $\kappa = |V_m|$.

\vspace{0.1cm}

\noindent Then in the second step: not all  $\kappa$-partitions of discrete boxes are to be considered.  \textcolor{blue}{\textbf{Only these}} partitions of $m$-dimensional box $V_{m,n}$  are admitted  for which all sub-boxes \textbf{are of the form} $V_m$  (\textcolor{blue}{\textbf{self-similarity}}).

\begin{defn}
Let $V_{m,n}$ be a $m$-dimensional $F$-box. Then any $\kappa$-partition into sub-boxes of the form $V_m$ is called tiling of $V_{m,n}$.
\end{defn}

\noindent It was shown in [13] that just the admissibility condition (see also [15]) is not sufficient for the existence a tiling for any given $m$-dimensional box $V_{m,n}$. Kwa\'sniewski in his papers [11,8]  posed the following problem called \emph{Cobweb Tiling Problem}, which was a starting point of the present author  research with  results being reported in the presents note.

\vspace{0.2cm}

\noindent \textbf{Tiling problem}\\
Suppose that $F$ is an admissible sequence. Under which conditions any $F$-box $V_{m,n}$ designated by sequence $F$ has a tiling? Find effective characterizations and/or find an algorithm to produce these tilings.

\vspace{0.1cm}

\noindent In order to recognize the present state of such cobwebs tiling investigation see Dziemia\'nczuk's papers [13.15.48] and contact [\textbf{17}] as well as  Internet Gian Carlo Rota Polish Seminar   \emph{http://ii.uwb.edu.pl/akk/sem/sem\_rota.htm}  (subjects \textbf{3}  and \textbf{5}).

%%%%%%%%%%%%%%%%%%%%%%%%%%%%%%%%%%%%%%%%%%%%%%%%%%%%%%%%%%%%%%%%%%%%%%%%%%%%%%%%%

\vspace{0.2cm}

\noindent \textbf{2.2. Influence posets model interpretation of  graded posets.}

%%%%%%%%%%%%%%%%%%%%%%%%%%%%%%%%%%%%%%%%%%%%%%%%%%%%%%%%%%%%%%%%%%%%%%%%%%%%%%%%%
%%%%%%%%%%%%%%%%%%%%%%%%%%%%%%%%%%%%%%%%%%%%%%%%%%%%%%%%%%%%%%%%%%%%%%%%%%%%%%%%%
\vspace{0.1cm}

\noindent Here we present  Agnarsson picture [53] of \textbf{influence poset} among the agents adopted to the notation and nomenclature of this note.

\noindent Namely, consider the number $\textcolor{red}{\textbf{s}}$ of "\textcolor{red}{\textbf{s}}ecret agents" $A_1, . . . ,A_s$  (or "Atoms" or genes or viruses or galaxes or terrorists [54]) as being  monitored over discrete times $t = 0, 1, . . . ,n$, $n \in N \cup\left\{\infty\right\}$. One may introduce a graded poset structure   consisting of the $s(n + 1)$ vertices  $A_k(t)$ in Hasse digraph , where a directed edge from   $A_k(t)$  up  to $A_k(t+1)$  is present if and only if an agent  $A_p$ has influenced an agent $A_q$ during the time interval from  $t$ to   $t + 1$. 

\vspace{0.1cm}

\noindent The resulting poset is  called the \textcolor{blue}{\textbf{influence poset}} among the agents [53] or terrorists [54].

\vspace{0.1cm}

\noindent What we have now are $n+1$ levels of the influence poset $P$ with single levels $\Phi_t = {A_1(t), . . . ,A_s(t)}$ - one for each  time $t = 0, 1, . . . ,n$; (up direction - time direction).

\noindent If we admit the number $\textcolor{red}{\textbf{s}}$ of "\textcolor{red}{s}pies"  to be varying with time according to the schedule:  $\left|\Phi_r\right| = r_F$, then we arrive at  an $F$-denominated graded poset. If in addition we assume the total influence (invigilation?)  i.e. each agent  $A_p$  influences every agent $A_q$ during the time interval from  $r$ to $r + 1$, then we are at home of cobweb posets $\Pi_n$ ,  $\Pi_\infty \equiv \Pi$.

%%%%%%%%%%%%%%%%%%%%%%%%%%%%%%%%%%%%%%%%%%%%%%%%%%%%%%%%%%%%%%%%%%%%%%%%%%%%%%%%%%%%%%%%%%%%%%%

%%%%%%%%%%%%%%%%%%%%%%%%%%%%%%%%%%%%%%%%%%%%%%%%%%%%%%%%%%%%%%%%%%%%%%%%%%%%%%%%%%%%%%%%%%%%%%%

%%%%%%%%%%%%%%%%%%%%%%%%%%%%%%%%%%%%%%%%%%%%%%%%%%%%%%%%%%%%%%%%%%%%%%%%%%%%%%%%%%%%%%%%%%%%%%%

%%%%%%%%%%%%%%%%%%%%%%%%%%%%%%%%%%%%%%%%%%%%%%%%%%%%%%%%%%%%%%%%%%%%%%%%%%%%%%%%%%%%%%%%%%%%%%%

%%%%%%%%%%%%%%%%%%%%%%%%%%%%%%%%%%%%%%%%%%%%%%%%%%%%%%%%%%%%%%%%%%%%%%%%%%%%%%%%%%%%%%%%%%%%%%%

%%%%%%%%%%%%%%%%%%%%%%%%%%%%%%%%%%%%%%%%%%%%%%%%%%%%%%%%%%%%%%%%%%%%%%%%%%%%%%%%%%%%%%%%%%%%%%%

\section{Zeta and inverse zeta functions and matrices formulas. Then Whitney numbers.}

\vspace{0.1cm}

\noindent \textbf{3.1. Zeta and inverse zeta functions in Knuth notation - natural labeling choice.}
\vspace{0.1cm}

\noindent \textbf{Knuth notation}. In the  wise Knuth's "notationlogy" note [24]   one finds among others the notation just for the purpose here: 

$$ [ s ] =\left\{ \begin{array}{cl} 1&if \ s \ is \  true ,\\0&otherwise.

\end{array} \right. $$

\vspace{0.1cm} 

\noindent Consequently for any set or class

$$[x=y] \equiv \delta (x,y). $$

\vspace{0.1cm} 

\noindent Consequently for $x,y,k,s,n \in N \cup \left\{\textcolor{blue}{\textbf{0}}\right\}$

$$[x<y] \equiv \sum_{k\geq \textcolor{red}{\textbf{1}}} \delta (x+k,y),$$

$$[x\textcolor{red}{\leq} y] \equiv \sum_{k\geq \textcolor{red}{\textbf{0}}} \delta (x+k,y).$$

\noindent Using this makes the authors expressions [3,2,4] of the $\zeta$ and $\mu = \zeta^{-1}$ in terms of  $\delta$ transparent and handy  if rewritten in Donald Ervin Knuth's  notation [24]. Namely:

\vspace{0.1cm} 
\noindent \textbf{Zeta function $\zeta \in I(\Pi,R)$ formulas.} [3,2]

$$\zeta(x,y) =\zeta_{\textcolor{green}{\textbf{1}}}(x,y) -\zeta_{0}(x,y)$$

$$\zeta_{\textcolor{green}{1}}(x,y)= [x \leq y]$$ 

$$\zeta_{0}(x,y)=\sum_{s \geq 1}\sum_{k \geq \textcolor{red}{\textbf{1}}}[x=k+s_F][1 \leq y \leq s_F + (s-1)_F\! - 1 ].$$ 

$$\zeta_{0}(x,y)=\sum_{s \geq 1}[x>s_F][1 \leq y \leq s_F + (s-1)_F\! - 1 ],$$ 
where, let us  recall natural labeling choice: $x,y,k,s \in N \cup \left\{\textcolor{blue}{\textbf{0}}\right\}$.
\noindent Note, that for \textcolor{red}{\textbf{$F$ = Fibonacci}} this still more simplifies   as then  

$$s_F + (s-1)_F-1 = (s+1)_F.$$
Finally - as derived in [3,2] - we have  for  \textbf{(}$x,y,k,s,n \in N \cup \left\{\textcolor{blue}{\textbf{0}}\right\}$

$$ \zeta(x,y) =  [x\leq y] - [x<y] \sum_{n\geq \textcolor{blue}{\textbf{0}}} [(x> S(n)] [y \leq S(n+1)].$$
where 
$$ S(n) = \sum_{k\geq 1}^n k_F $$
\textbf{The M{\"{o}}bius $\mu$ i.e. inverse zeta function $\zeta^{-1} \in I(\Pi,R)$ formula} [3,2].

\vspace{0.1cm}

\begin{figure}[ht]
\begin{center}
	\includegraphics[width=70mm]{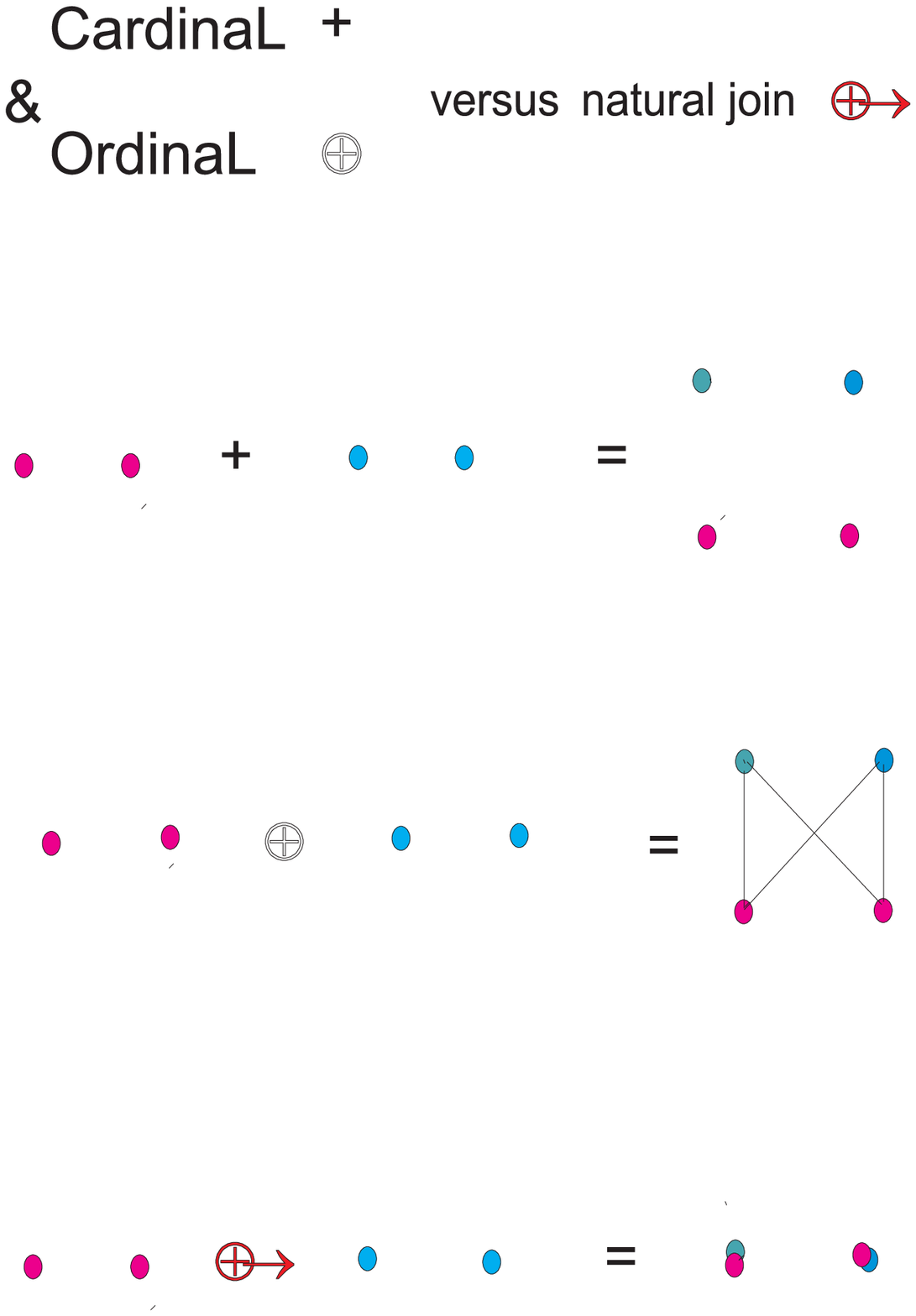}
	\caption{Cardinal and ordinal sums versus natural join.} \label{fig:representation}
\end{center}
\end{figure}

\noindent Now we shall quote  the M{\"{o}}bius $\mu \in I(\Pi,R)$ formula for $F$-denominated cobweb posets from [3,2] to where we refer for the derivation and discussion of this formula history. In its presentation we shall use the  plane grid coordinates description of cobweb posets Hasse diagrams  introduced in 2003 by the author of [18,19,20]  for the definition of now KoDAGs  and then  consequently used by the author of [25-29].  Here comes this formula right after the notation of  plane grid coordinates.\\
Namely, we label vertices of poset $P$ Hasse diagram  by pairs of coordinates: $\langle s,t \rangle \in {\bf N}\times {\bf N_{0}}$, where  $t$ coordinate is the level label while $N_{0} = N\cup\left\{0\right\}$. Then 
$$\Phi_{t}=\left\{\langle s,t \rangle ,\;\;1\leq s \leq t_F\right\},\;\;\;t\in{\bf N}\cup\{0\}.$$
Consequently  Hasse digraph of the poset  $P$ is a labeled digraph $D_P=\left(V,E\right)$ where

$$V=\bigcup_{t\geq 0}\Phi_{t},\;\;\;E=\left\{\langle \,\langle s,t\rangle,\langle r,t+1
\rangle\,\rangle\right\},\;\;1\leq s\leq t_F,\;\;1\leq r\leq (t+1)_F.$$

\noindent The partial order relation on $P$ is then described as follows: for  $x=\langle s,t\rangle, y=\langle u,v\rangle\;\;\;x,y\in P $ 
$$ ( x \leq_{P} y) \Longleftrightarrow  [(t<v)\vee (t=v \wedge s=u)].$$

%%%%%%%%%%%%%%%%%%%%%%%%%%%%%%%%%%%%%%%%%%%%%%%%%%%%%%%%%%%%%%%%%%%%%%%%%

\noindent Observe now  sine qua non  conditions $(*)$ - which stem from the definition  of  $\Phi_t$ level  i.e.  $(*)$ means:  $1 \leq s \leq t_F$, $1 \leq u \leq v_F$ while $t,v \in N_0$.  Then  with these sine qua non  conditions being \textbf{implemented in the formula} we get

$$\mu(x,y) = \mu (\left\langle s,t\right\rangle,\left\langle u,v\right\rangle) = [(s=u]\: [t=v]  - [t+1=v]+$$  
$$ +  [t+1<v]\:[1 \leq s \leq t_F ] \:[1 \leq u \leq v_F ]
(-1)^k\prod_{i= t+1}^{v-1}(i_F-1).$$
The proof of this formula relies on and consists of the derivation of is version in  \textbf{1.4. Cobweb posets' M\"{o}bius function} 
where it is the obvious statement for  $s>r , \ [x,z] = x\oplus \Phi_{r+1}\oplus...\oplus \Phi_{s-1}\oplus z$   ( conditions $(*)$ not  implemented in it) i.e. 

  $$\mu(x,y) = (-1)^{s-r}\prod_{k=r+1}^{s-1}[k_F-1] $$
($\left|\Phi_k\right|= k_F$,  $x,y \in P$) - the statement to be just rewritten in plane  grid coordinates  and with these $(*)$ sine qua non  \textbf{conditions} being implemented in the formula.

%%%%%%%%%%%%%%%%%%%%%%%%%%%%%%%%%%%%%%%%%%%%%%%%%%%%%%%%%%%%%%%%%%%%%%%%%

%%%%%%%%%%%%%%%%%%%%%%%%%%%%%%%%%%%%%%%%%%%%%%%%%%%%%%%%%%%%%%%%%%%%%%%%%

%%%%%%%%%%%%%%%%%%%%%%%%%%%%%%%%%%%%%%%%%%%%%%%%%%%%%%%%%%%%%%%%%%%%%%%%%

\vspace{0.2cm}

\noindent \textbf{3.2. Zeta and inverse zeta \textcolor{blue}{matrices}  - natural labeling choice.}
\vspace{0.1cm}

\noindent The matrix elements of $\zeta(x,y)$ matrix for Fibonacci cobweb poset were given in 2003 ([19,23] Kwa\'sniewski) using   $x,y \in N \cup \left\{\textcolor{blue}{\textbf{0}}\right\}$ labels of  vertices in their natural linear extension order i.e.  applying the natural labeling (see [55] - see Fig.11),  where - recall: \textit{The poset $P$ is naturally labeled if $x_i <x_j$ in $P$ implies that $i < j$}.  For that to do proceed as follows   (see Fig.11). 

\vspace{0.2cm}

\noindent \textbf{1.} set  \textcolor{blue}{\textbf{$k=0$}},\\ 
\textbf{2.} then label subsequent vertices - from the left to the right - along the level  $k$,\\
\textbf{3.} repeat  2.  for $k \rightarrow  k+1$  until  $k=n+1$ ;  $n  \in N\cup\left\{\infty\right\}.$ 
\noindent The labeling $\omega : \Pi \longrightarrow {0,1,..., n}$, $n \in N\cup\left\{\infty \right\}$  we thus arrive at  is  natural [55] i.e. $\omega$  is order preserving.

\vspace{0.2cm}
\noindent As the result we obtain for example  the $\zeta$ matrix for Fibonacci sequence as presented by the the Example.1 with \textit{La Scala di Fibonacci} in [3,4,2] and  dating back to 2003 [19,23]. For many other examples see [3,2]. For general case zeta matrix formula i.e. for  any $F$-graded poset see [4,3,2] or see Theorem 4 and make proper extension yourself. Here we quote the explicit  expression of \textbf{cobweb posets'} zeta matrix $\zeta_F$  for arbitrary natural numbers valued $F$- sequence  derived in [4] due to more than  mnemonic  efficiency  of the up-side-down notation being applied. With this notation inspired by Gauss numbers $k_q$ and replacing  $k$ - natural numbers with "$k_F$ numbers" = elements of  the  $F $-sequence one gets the theorem.

\vspace{0.3cm}

\noindent \textbf{Theorem 3} (Kwa\'sniewski)

$$
	\zeta_F = exp_\copyright[\mathbf{A}_F] \equiv (1 - \mathbf{A}_F)^{-1\copyright} \equiv I_{\infty\times\infty} + \mathbf{A}_F + \mathbf{A}_F^{\copyright 2} + ... =
$$
$$
	= \left[\begin{array}{lllll}
	\textcolor{red}{\textbf{I}}_{1_F\times 1_F} & I(1_F\times\infty) \\
	O_{2_F\times 1_F} & \textcolor{red}{\textbf{I}}_{2_F\times 2_F} & I(2_F\times\infty) \\
	O_{3_F\times 1_F} & O_{3_F\times 2_F} & \textcolor{red}{\textbf{I}}_{3_F\times 3_F} & I(3_F\times\infty) \\
	O_{4_F\times 1_F} & O_{4_F\times 2_F} & O_{4_F\times 3_F} & \textcolor{red}{\textbf{I}}_{4_F\times 4_F} & I(4_F\times\infty) \\
	... & etc & ... & and\ so\ on & ...
	\end{array}\right]
$$
where  $I (s\times k)$  stays for $(s\times k)$  matrix  of  ones  i.e.  $[ I (s\times k) ]_{ij} = 1$;  $1 \leq i \leq  s,  1\leq j  \leq k.$  and  $n \in N \cup \{\infty\}$ and where

$$
	\mathbf{A}_F = \left[\begin{array}{llllll}
	0_{1_F\times 1_F} & I(1_F \times 2_F) & 0_{1_F \times \infty} \\
	0_{2_F\times 1_F} & 0_{2_F\times 2_F} & I(2_F \times 3_F) & 0_{2_F \times \infty} \\
	0_{3_F\times 1_F} & 0_{3_F\times 2_F} & 0_{3_F\times 3_F} & I(3_F \times 4_F) & 0_{3_F \times \infty} \\
	0_{4_F\times 1_F} & 0_{4_F\times 2_F} & 0_{4_F\times 3_F} & 0_{4_F\times 4_F} & I(4_F \times 5_F) & 0_{4_F \times \infty} \\
	... & etc & ... & and\ so\ on & ...
	\end{array}\right]
$$
\noindent In the $\zeta_F $ formula from [6,4]  $\copyright$ denotes the Boolean  product, hence - exactly this product is meant while  Boolean powers enter formulas. We readily recognize from its block structure that $F$-La Scala descending and descending far away own to infinity -  is formed by \textcolor{red}{\textbf{upper zeros}} of block-diagonal matrices $\textcolor{red}{\textbf{I}}_{k_F\times k_F}$.

\vspace{0.2cm}

\noindent \textcolor{blue}{\textbf{Extension of the Theorem 3.}} In the Theorem 3  $I(k_F \times (k+1)_F)$  denotes  $k_F \times (k+1)_F$ matrix  of all entries equal to one. \textbf{For any} $F$-\textbf{denominated poset} \textcolor{red}{\textbf{replace}} $I(k_F \times (k+1)_F)$ \textcolor{red}{\textbf{by}} $B(k_F \times (k+1)_F)$ obtained from  $I(k_F \times (k+1)_F)$ via replacing adequately (in accordance with Hasse  digraph) corresponding \textit{ones} by \textbf{zeros}.

\vspace{0.2cm}

\noindent The next formula - for inverse zeta function $\zeta^{-1} = \mu$  as represented by matrix in natural labeling of $F$-graded poset - was derived in [3,2]  and is quoted below  as the Theorem 4.

\vspace{0.4cm}

\begin{flushright}

1 - - - - - - - - - - - - - - - - - - - - - - - - - - - - - - - - - - - - - - - - - - - - - - - - - - -\\
1 - - - - - - - - - - - - - - - - - - - - - - - - - - - - - - - - - - - - - - - - - - - - - - - - - -\\
1\textbf{ \textcolor{red}{0}} - - - - - - - - - - - - - - - - - - - - - - - - - - - - - - - - - - - - - - - - - - - - - - -\\
1 - - - - - - - - - - - - - - - - - - - - - - - - - - - - - - - - - - - - - - - - - - - - - - -\\
1 \textbf{\textcolor{red}{0 0}} - - - - - - - - - - - - - - - - - - - - - - - - - - - - - - - - - - - - - - - - - - -\\
1 \textbf{\textcolor{red}{0}} - - - - - - - - - - - - - - - - - - - - - - - - - - - - - - - - - - - - - - - - - - -\\
1 - - - - - - - - - - - - - - - - - - - - - - - - - - - - - - - - - - - - - - - - - - -\\
$F_{5}-1\; \textbf{\textcolor{red}{0}}'s \;\; $1 \textbf{\textcolor{red}{0 0 0 0}} - - - - - - - - - - - - - - - - - - - - -  - - - - - - - - - - - - - - -\\
1 \textbf{\textcolor{red}{0 0 0}} - - - - - - - - - - - - - - - - - - - - - - - - - - - - - - - - - - - -\\
1 \textbf{\textcolor{red}{0 0}} - - - - - - - - - - - - - - - - - - - - - - - - - - - - - - - - - - - -\\
1 \textbf{\textcolor{red}{0}} - - - - - - - - - - - - - - - - - - - - - - - - - - - - - - - - - - - -\\
1 - - - - - - - - - - - - - - - - - - - - - - - - - - - - - - - - - - - -\\
$F_{6}-1\;zeros \quad \quad \quad \quad $1 \textbf{\textcolor{red}{0 0 0 0 0 0 0}} - - - - - - - - - - - - - - - - - - - - - - - - - -\\
1 \textbf{\textcolor{red}{0 0 0 0 0 0}} - - - - - - - - - - - - - - - - - - - - - - - - -\\
1 \textbf{\textcolor{red}{0 0 0 0 0}} - - - - - - - - - - - - - - - - - - - - - - - - -\\
1 \textbf{\textcolor{red}{0 0 0 0}} - - - - - - - - - - - - - - - - - - - - - - - - -\\
1 \textbf{\textcolor{red}{0 0 0}} - - - - - - - - - - - - - - - - - - - - - - - - -\\
1 \textbf{\textcolor{red}{0 0}} - - - - - - - - - - - - - - - - - - - - - - - - -\\
1 \textbf{\textcolor{red}{0}} - - - - - - - - - - - - - - - - - - - - - - - - -\\
1 - - - - - - - - - - - - - - - - - - - - - - - - -\\
$F_{7}-1\;zeros \quad \quad \quad \quad \quad \quad \quad \quad
\quad \quad   \;$1 \textbf{\textcolor{red}{\textcolor{red}{0 0 0 0 0 0 0 0 0 0 0 0}}} - - - - - - - \\
1\textbf{ \textcolor{red}{0 0 0 0 0 0 0 0 0 0 0}} - - - - - - -\\
1\textbf{ \textcolor{red}{0 0 0 0 0 0 0 0 0 0}} - - - - - - -\\
1 \textbf{\textcolor{red}{0 0 0 0 0 0 0 0 0}} - - - - - - -\\
1\textbf{ \textcolor{red}{0 0 0 0 0 0 0 0}} - - - - - - - \\
1\textbf{ \textcolor{red}{0 0 0 0 0 0 0}} - - - - - - - \\
1 \textbf{\textcolor{red}{0 0 0 0 0 0}} - - - - - - - \\
1 \textbf{\textcolor{red}{0 0 0 0 0}} - - - - - - - \\
1 \textbf{\textcolor{red}{0 0 0 0}} - - - - - - - \\
1 \textbf{\textcolor{red}{0 0 0}} - - - - - - - \\
1 \textbf{\textcolor{red}{0 0}} - - - - - - -  \\
1
\textbf{\textcolor{red}{0}} - - - - - - - \\
1 - - - - - - - \\
$F_{8}-1\;zeros \quad \quad \quad \quad \quad \quad \quad \quad
\quad \quad \quad \quad \quad \quad \quad \quad \quad \quad \quad \quad \;\,$.........................\\
{\em and so on}
\end{flushright}

\vspace{0.1cm}
\begin{center}
\noindent \textbf{ Example.1  La Scala di Fibonacci . The  staircase structure  of  incidence
matrix $\zeta_F$}  for \textbf{$F$=Fibonacci sequence} \end{center}

\vspace{0.1cm}

\noindent \textbf{Theorem 4} (Kwa\'sniewski)

\vspace{0.1cm}

\noindent Let $F$ be \textbf{any} natural numbers valued sequence. Then \textcolor{red}{\textbf{for arbitrary}} $F$-denominated graded poset (cobweb posets included)    

$$ C(\mu_F)_{r,s}= c_{r,s} = [r=s] + K_s(r_F)(-1)^{s-r} =[r=s] + [s>r] (-1)^{s-r}[(r+1)_F - 1]^{\overline{s-r}} $$ 
i.e.

$$C(\mu_F)_{r,s}= c_{r,s} = [s\geq r] (-1)^{s-r}[(r+1)_F - 1]^{\overline{s-r}}, $$
i.e. 
$$C(\mu_F)_{r,s}= c_{r,s} = [s\geq r] (-1)^{s-r}K_s(r_F), $$

\noindent with matrix elements from $N$ or the ring  $R$= $2^{\left\{1 \right\}}$ , $Z_2=\left\{0,1\right\}$, $Z$ etc.\\
i.e. for cobweb posets

$$\mu = \left( \delta_{r,s}I_{r_F\times r_F} + (-1)^{s-r}K_s(r_F)I(r_F \times s_F)   \right)$$ 
i.e.

$$
	\mu  = \left[\begin{array}{llllll}
	I_{1_F\times 1_F} & c_{1,2}I(1_F \times 2_F) & c_{1,3}I(1_F \times 3_F) & c_{1,4}I(1_F \times 4_F) & c_{1,5}I(1_F \times 5_F) & c_{1,6}I(1_F \times 6_F)\\
	0_{2_F\times 1_F} & I_{2_F\times 2_F} & c_{2,3}I(2_F \times 3_F) & c_{2,4}I(2_F \times 4_F) & c_{2,5}I(2_F \times 5_F) & c_{2,6}I(2_F \times 6_F)\\
	0_{3_F\times 1_F} & 0_{3_F\times 2_F} & I_{3_F\times 3_F} & c_{3,4}I(3_F \times 4_F) & c_{3,5}I(3_F \times 5_F) & c_{3,6}I(3_F \times 6_F)\\
	0_{4_F\times 1_F} & 0_{4_F\times 2_F} & 0_{4_F\times 3_F} & I_{4_F\times 4_F} & c_{4,5}I(4_F \times 5_F) & c_{4,6}I(4_F \times 6_F)\\
	... & etc & ... & and\ so\ on & ...
	\end{array}\right]
$$
where  $I(k_F \times (k+1)_F)$  denotes (recall)  $k_F \times (k+1)_F$ matrix  of all entries equal to one. \textbf{For any} $F$-\textbf{denominated poset} \textcolor{red}{\textbf{replace}} $I(k_F \times (k+1)_F)$ \textcolor{red}{\textbf{by}} $B(k_F \times (k+1)_F)$ obtained from  $I(k_F \times (k+1)_F)$ via replacing adequately (in accordance with Hasse  digraph) corresponding \textit{ones} by \textbf{zeros}.

\noindent The functions' upper factorial from [2] used above  (valid whenever defined for corresponding functions $f$ of the natural number argument or of an argument from any chosen ring  -  Z,R,etc.) is defined as follows

$$
f(r_F)^{\overline{k}} = f(r_F)f([r+1]_F)...f([r + k -1]_F), \ n^{\overline{0}} \equiv 1,\: n \in N \cup \ \left\{ \textit{\textbf{0}} \right\},Z,R,etc. , 
$$

\vspace{0.1cm}

\noindent \textbf{3.3. Whitney numbers and characteristic polynomials for cobweb posets.}[2]\\

\vspace{0.1cm}

\noindent Let us remind the notation: $r_F = \left|\Phi_r\right|$. Let us then recall (see [56]) definitions of  Whitney numbers of the first kind $w_r(P)$ and Whitney numbers of second kind $W_r(P)$ where $P$ is any given graded poset with bounded independence sets - i.e.  $\left|\Phi_r\right| \in N$ for $r \in N \cup \left\{0\right\}$  and we assume that $\left|\Phi_0\right|= 1 = 0_F$  hence  $\Phi_0 = \left\{0\right\}$; $0 \in P$ stays for minimal element of the poset $P$. Now here are these definitions. 

$$w_{r}(P)=\sum_{x\in P,\,r(x)=r}\mu(0,x),$$
$$W_r(P)=\sum_{x\in P,\,r(x)=r}1=|\{x\in P\,:\,r(x)=r\}|.$$

\vspace{0.2cm}

\noindent It is obvious just by notation  that 

$$W_r(P)= r_F .$$
Of course in the general case of finite posets $P_n$ the number $\left|\Phi_k\right|$ might depend on $n$ and then we end up with an array $\left(W_k(P_n)\right)$ of Whitney numbers as it is the case with unimodal binomials or unimodal Gaussian binomials for example  or Stirling-like numbers as in [45] or...- see further classical examples in [56,57,58).

\vspace{0.2cm}

\noindent According to Theorem 4 we have for \textcolor{blue}{\textbf{cobweb posets}} i.e. for $\Pi$'s

$$C(\mu_F)_{0,s}= c_{0,s} = [s\geq 0] (-1)^s(1_F - 1)^{\overline s}, $$
or equivalently - as the values of $\mu(x,y)$ depend only on the rank of its arguments  
$$\mu_F(0,x)=  [r(x)\geq 0] (-1)^{r(x)}(1_F - 1)^{\overline {r(x)}}, $$
or equivalently (compare all this with (3) in [29])
$$\mu_F(0,x)=  [x=0] - [r(x)=1](1_F - 1) +  [r(x)>1](-1)^{r(x)}\prod_{k=1}^{r(x)-1}(k_F-1), $$
or equivalently  just 
$$\mu_F(0,x)=  [r(x)\geq 0] (-1)^{r(x)}K_r(0_F). $$

\vspace{0.1cm}

\noindent Consequently - as the values of $\mu(x,y)$ depend only on the rank of its arguments - the Whitney numbers of the first kind for the denominated by $F$ 
\textcolor{blue}{\textbf{cobweb poset}} $\Pi$ may be calculated along the formula

  $$w_r(\Pi) =  \sum_{\{x\in\Pi\,:\,r(x)=r\}}\mu_F(0,x)= r_F \cdot\mu_F(0,x)$$
 i.e.
   
   $$w_r(\Pi) = r_F\cdot(-1)^r K_r(0_F).$$
Naturally $w_0(\Pi)=1$. Compare  the above with (4) in [29].  

\vspace{0.2cm}

\noindent Of course in the \textit{general case} of finite posets $P_n = \bigcup_{k=0}^n \Phi_k(n) $ the number $\left|\Phi_k(n)\right|$ might depend on $n$ and then we end up with an array $\left(w_k(P_n)\right)$ of Whitney numbers of the first kind as it is the case with binomials or Gaussian binomials for example - (see further classical examples in [56,57,58]).

\vspace{0.2cm}

\noindent To this end - consequently - let us consider characteristic polynomials  $\chi_{P_n}(t)$, $n\geq 0$ 
defined as ([56,57], [29])

$$\chi_{P_n}(t)=\sum_{x\in P_n}\mu(0,x)t^{n-r(x)}=\sum_{k = 0}^n w_k(P_n)t^{n-k}.$$
The formula for characteristic polynomials - here for \textbf{specific} $F$-denominated  finite cobweb sub-posets    $P_n = \oplus_{k=0}^n \Phi_k$   (i.e. $\left|\Phi_k(n)\right|$ does not depend on $n$)- namely - this  formula for characteristic polynomials  obviously is of the form 

$$ \chi_{P_n}(t)=\sum_{k=0}^n (-1)^k \cdot k_F \cdot x^{n-k} \cdot K_k(0_F) $$
or equivalently  (compare with Theorem 3.1 in [29])

$$ \chi_{P_n}(t)=x^n - x^{n-1}1_F(1_F-1)  +  \sum_{k=2}^n (-1)^k \cdot k_F \cdot x^{n-k} \cdot \prod_{r=1}^{k-1}(r_F-1).$$

\vspace{0.1cm}

\begin{figure}[ht]
\begin{center}
	\includegraphics[width=70mm]{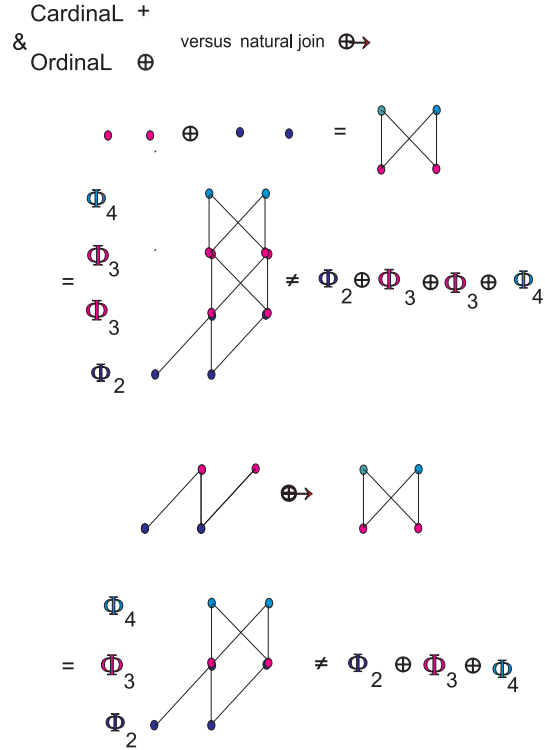}
	\caption{Ordinal sum versus natural join.} \label{fig:representation}
\end{center}
\end{figure}

\begin{figure}[ht]
\begin{center}
	\includegraphics[width=70mm]{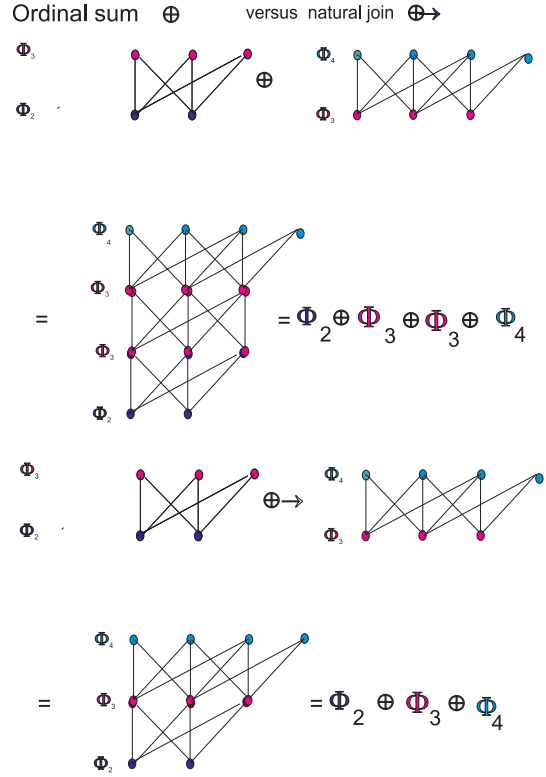}
	\caption{Ordinal sum versus natural join - again.} \label{fig:representation}
\end{center}
\end{figure}

\vspace{0.1cm}

\noindent \textbf{Recapitulation  of \textcolor{red}{natural join} ideology.}

\vspace{0.2cm}

\noindent \textbf{Recall} that  both $\leq$ partial order relations and  $\prec\cdot$ cover relations are \textcolor{red}{\textbf{natural join}} of their bipartite \textcolor{blue}{\textbf{relations}} correspondent \textcolor{blue}{\textbf{chains}}, and this is exactly the reason and the very source of the Theorem 4  validity and shape. Note also that all information on structure of any  poset $P$ is 
coded by the  $\zeta$ matrix - a characteristic function of $\leq \in \  P = \left\langle \Phi,\leq\right\rangle$. In short: $\zeta$ and equivalently $\mu = \zeta^{-1}$ are the Incidence algebra of $P$ 
coding elements. In brief - the following identifications are self-evident:

$$\left\langle \Phi,\mu_F \right\rangle   \equiv \left\langle \Phi,\zeta_F \right\rangle  \equiv  \left\langle \Phi,\leq \right\rangle \equiv  \left\langle \Phi,\textcolor{blue}{\textbf{C}}(\mu_F) \right\rangle \equiv  \left\langle \Phi,\prec\cdot \right\rangle.$$
In case of  $F$-graded posets all of them i.e.  $\mu_F, \zeta_F, \prec\cdot, \leq $ are of \textcolor{red}{\textbf{natural join}} of their bipartite counterparts origin [6,5,4,3,2]. As a matter of illustration we quote the definition (6 in [6]) of cobweb poset from [6]. 

\begin{defn} [cobweb poset]
Let $D = (\Phi, \prec\!\!\cdot )$  be a  transitive irreducible  digraph. Let  $n \in N \cup \{\infty\}$.  Let $D$ be a natural join $D = \os_{k=0}^n B_k$   of   di-bicliques  $B_k = (\Phi_k \cup \Phi_{k+1}, \Phi_k\times\Phi_{k+1} ) , n \in N \cup \{\infty\}$.    Hence the digraph  $D = (\Phi,\prec\!\!\cdot  )$    is graded. The poset $\Pi (D)$ associated to this graded digraph $D = (\Phi,\prec\!\!\cdot )$    is called a cobweb poset.
\end{defn}

%%%%%%%%%%%%%%%%%%%%%%%%%%%%%%%%%%%%%%%%%%%%%%%%%%%%%%%%%%%%%%%%%55

\vspace{0.1cm}

\section{$\os$ versus $+$  and $\oplus$. Summary of elementary properties and related statements}

\vspace{0.1cm}

\noindent \textbf{4.1 \textcolor{blue}{The Principal properties of the natural join operation}.}

\vspace{0.1cm}
\noindent \textbf{4.1.1.} Particular cases (see 1.2.5.). Whether posets must be or need not to be finite - becomes clear from the context.

\vspace{0.1cm}
\noindent Convention:  $P\oplus \emptyset) \oplus (\emptyset \oplus P)= P$.

\vspace{0.1cm} 

\noindent  Let \textbf{nonempty posets} $\left\langle P,\leq_P \right\rangle$ and $\left\langle Q,\leq_Q \right\rangle$ do satisfy the \textcolor{blue}{\textbf{Natural Join Conformity Condition}} i.e. are mutually \textbf{not disjoint} or constitute posets with partial orders $\leq_P $  and $\leq_Q$ being equivalent on corresponding in $P$ and $Q$ isomorphic sub-posets  which we shall consider identical by convention. Let then $P =P_1 \cup P_2$  and  $Q = P_2 \cup P_3$. The $P\os Q $ is  defined via Definition 3. Here we recall  some special cases to be used right next.

\vspace{0.1cm}
\noindent Let $P =P_1 + P_2$  and  $Q = P_2 + P_3$ then we define $\os$ via identity

\begin{center}
$P \os Q  \equiv P_1 + P_2 + P_3 $
\end{center}

\vspace{0.1cm}
\noindent Let $P =P_1 + P_2$  and  $Q = P_2 \oplus P_3$ then we define $\os$ via identity
\vspace{0.1cm}
\begin{center}
$P \os Q  \equiv P_1 + P_2 \oplus P_3 $
\end{center}

\vspace{0.1cm}

\noindent  Let  $\left\langle P_1, \leq_1 \right\rangle$ , $\left\langle P_2, \leq_2\right\rangle$,  $\left\langle P_3, \leq_3\right\rangle$  be posets. Let $P =P_1 \cup P_2$  and  $Q = P_2\cup P_3$ be disjoint sums i.e. represent corresponding two block partitions of $P$  and  $Q$.  Let $x \leq_1y \equiv  x \leq_2y$  for  $x,y \in P_2$.  Then

\vspace{0.1cm}
\noindent Let $P =P_1 \oplus P_2$  and  $Q = P_2 \oplus P_3$ then we define $\os$ via identity

\begin{center}

$P \os Q  \equiv \Pi_1 \oplus \Pi_2 \oplus \Pi_3 $
\end{center}

\begin{center}
$  \left\langle P \os Q \leq_{\vspace{0.1cm}}\right\rangle  \equiv \left\langle P_1 \cup P_2 \cup P_3,\leq_{\vspace{0.1cm}}\right\rangle $
\end{center}
where $\leq_{P \os Q }$ is defined via Definition 3 in 1.2.5.

%%%%%%%%%%%%%%%%%%%%%%%%%%%%%%%%%%%%%%%%%%%%%%%%%%%%%%%%%%%%%%%%%%%%%%%%%%%%%%

%%%%%%%%%%%%%%%%%%%%%%%%%%%%%%%%%%%%%%%%%%%%%%%%%%%%%%%%%%%%%%%%%%%%%%%%%%%%%%

\vspace{0.2cm}

\noindent \textbf{4.1.2.} Recall now the elementary properties of cardinal $+$  and ordinal $\oplus$ sums  [59,1]. Let $X.Y,Z$ be posets. \textcolor{red}{\textbf{Let}} us define for   
$\left\langle X,\prec \right\rangle$ the  $\left\{X\right\}=$ number of pairs in covering relationship  $\prec \cdot$ which is induced by partial order relation $\prec$ = number of arrows in Hasse digraph $\left\langle X, \prec \cdot \right\rangle$ of the poset  $\left\langle X,\prec \right\rangle$. Then

\vspace{0.1cm}
\noindent  [\textbf{P.1}]  $X+Y = Y+X$ , [\textbf{P.2}]  $(X+Y)+ Z = X+(Y+Z)$ ,  [\textbf{P.3}]  $\left|X+Y\right| = \left|X\right|  +  \left|Y\right|$, [\textbf{P.4}]  $\left\{X+Y\right\} = \left\{X\right\}  +  \left\{Y\right\}$.\\
Let since now on $A[X]$ denotes the adjacency matrix of the poset $\left\langle X,\prec \right\rangle$. Then   the adjacency matrix $\mathbf{A}[X+Y]$  of a Hasse DAG  digraph $D(X+Y)$ is the direct sum of corresponding adjacency matrices $\mathbf{A}[X]$ and $\mathbf{A}[Y]$ \\  
\vspace{0.1cm}
\noindent [\textbf{P.5}]

$$
	\mathbf{A}[X+Y] = \left(
	\begin{array}{cc}
	\mathbf{X} & \mathbf{0} \\
\mathbf{0} & \mathbf{Y} \\
	\end{array}
	\right).
$$
\noindent \textcolor{red}{\textbf{Let}} us define now for posets $\left\langle X,\prec \right\rangle$ and $\left\langle Y,\prec* \right\rangle$  their  Cartesian product  $X \cdot Y $   
$\left\langle X \times Y,\leq  \right\rangle$ the  $\left\{X\right\}$  via  $ (x_1,y_1) \leq (x_2,y_2)$   iff  $x_1 \prec  x_2$  and  $y_1 \prec* y_2$. Then \\  
\noindent [\textbf{P.6}]  $X \cdot Y = Y \cdot X$ , [\textbf{P.7}]  $(X \cdot Y) \cdot Z = X \cdot (Y \cdot Z)$ and   [\textbf{P.8}]  $(X + Y) \cdot Z = X \cdot Z + Y \cdot Z$\\   
- (the same for lexicographic product,  which is also distributive with respect to $\oplus$ [59]) while [\textbf{P.9}] $(X \oplus Y) \cdot Z \neq  X \cdot Z \oplus Y $. For $\oplus$ other elementary properties are the  following  ones.  \textbf{P.10}  $X \oplus Y \neq Y \oplus X$  and  \textbf{P.11}  $(X \oplus Y) \oplus Z = X \oplus (Y \oplus Z)$.

\vspace{0.2cm}
\noindent The adjacency matrix $\mathbf{A}[X+Y]$  of a Hasse DAG  digraph $D(X\oplus Y)$ \textit{is not the direct sum} of corresponding adjacency matrices $\mathbf{A}[X]$ and $\mathbf{A}[Y]$.  

\vspace{0.1cm}
\noindent [\textbf{P.12}]

$$
	\mathbf{A}[X\oplus Y] = \left(
	\begin{array}{cc}
	\mathbf{X} & \mathbf{I(r \times s )} \\
\mathbf{0} & \mathbf{Y} \\
	\end{array}
	\right),
$$
where $r = \left|X\right|$, $s = \left|Y\right|$  and $I (r\times s)$  stays for  $(r\times s)$  matrix  of  ones  i.e.  $[ I (r\times s) ]_{ij} = 1$;  $1 \leq i \leq r,  1\leq j  \leq s$. 

\vspace{0.1cm}
\noindent  [\textbf{P.13}] $\left\{X\oplus Y\right\} = \left\{X\right\}  +  \left\{Y\right\}  +  M(X)\cdot m(Y)$ where  $M(X)$ is the number of maximal elements in $X$ and  $m(Y)$ is the number of minimal elements in $Y$ . Hence in general [\textbf{P.14}] $\left\{X\oplus Y\right\}\neq \left\{Y\oplus X\right\}$. Naturally [\textbf{P.15}]  $\left|X\otimes Y\right| = \left|X\right|  +  \left|Y\right|$.

\vspace{0.2cm}

\noindent \textbf{4.1.3.} In view of the above the following is obvious.

\vspace{0.1cm}
\noindent \textbf{Obvious.} The natural join operation $\os$ is associative and noncommutative - as the ordinal [linear] sum $\oplus$ is. 

\vspace{0.1cm}

\noindent [\textbf{P.16}]  $X \os Y \neq Y \os X$  ,    [\textbf{P.17}]  $(X \os Y) \os Z = X \os (Y \os Z)$ 

\vspace{0.1cm}
\noindent \textbf{Obvious.} Naturally \textbf{always}  [\textbf{P.18}] $P\os Q \neq P\oplus Q $ for non-empty  sets $P$  and $Q$ - see Figures 13,14,15.

\vspace{0.1cm}
\noindent [\textbf{P.19}] $(X \os Y) \cdot Z \neq  X \cdot Z \os Y \prec Z$,   [\textbf{P.17}] $\left\{X\os Y\right\} = \left\{X\right\} + \left\{Y\right\} - \left\{X\cap Y\right\}$ . Hence in general [\textbf{P.20}] $\left\{X\os Y\right\}\neq \left\{Y\os X\right\}$  if defined at all. Similarly  [\textbf{P.20}]  $\left|X\os Y\right| = \left|X\right|+ \left|Y\right| - \left|X\cap Y\right| $; see Figures 13,14,15.

\vspace{0.2cm}

\noindent [\textbf{P.21}] \textbf{Natural join of Hasse DAG diagrams adjacency matrices.} .

\vspace{0.1cm}
\noindent The adjacency matrix    $A[D(P)] \equiv A[D]$   of a DAG   $D(P)= \left\langle P,\prec\cdot = E\subseteq P\times P \right\rangle$  i.e. Hasse digraph of the poset $\left\langle P,\prec \right\rangle$ is given by

$$
	\mathbf{A}[D] = \left(
	\begin{array}{cc}
  0  &  P_{r,s} \\
\mathbf{zeros} & 0 \\
	\end{array}
	\right),
$$
where $r = 2,...,\left|P\right| \geq s > r$.

\vspace{0.1cm}
\noindent  Correspondingly the adjacency matrix   $A[D(P)] \equiv A[D]$ of a DAG such that  $D(P)= (V\cup W, E \subseteq  P\times P) )$ were $V\cup W = P$ and $V\cap W = \emptyset $ , is given in self-explanatory notation by

$$
	\mathbf{A}[D] = \left(
	\begin{array}{cc}
 \mathbf{A[D(V)]} &   A[D(V,W)] \\
\mathbf{zeros} &  \mathbf{A[D(W)]} \\
	\end{array}
	\right),
$$
where   $A[D(V,W)] = \left(A[D]_{r,s}\right)$  for $r = 1,...,\left|V\right|$ and  $s = 1,...,\left|W\right|$  and  $A[D]_{r,s} = [x_r \in V \;\wedge \:y_s \in W \;\wedge \;x_r \prec\cdot y_s ]  $ where Knuth notation [24] was used.

\begin {defn} {Natural join of Hasse DAG diagrams adjacency matrices.} 

\noindent Let \textcolor{blue}{\textbf{Natural Join Conformity Condition}} be  satisfied by $A[D(V\cup W)]$ and $A[D(W\cup U)]$ i.e. let
$$
	\mathbf{A}[D(V\cup W)] = \left(
	\begin{array}{cc}
 \mathbf{A[D(V)]} &   A[D(V,W)] \\
\mathbf{zeros} &  \mathbf{A[D(W)]} \\
	\end{array}
	\right),
$$
and
$$
	\mathbf{A}[D(W\cup U)] = \left(
	\begin{array}{cc}
 \mathbf{A[D(W)]} &   A[D(W,U)] \\
\mathbf{zeros} &  \mathbf{A[D(U)]} \\
	\end{array}
	\right),
$$ 
then  [\textbf{P.21}]
$$
	\mathbf{A}[D(V\cup W)]\os 	\mathbf{A}[D(W\cup U)] = \left(
	\begin{array}{ccc}
 \mathbf{A[D(V)]} &   A[D(V,W)] & 0\\
\mathbf{zeros} &  \mathbf{A[D(W)]} &   A[D(W,U)]\\
\mathbf{zeros} &  \mathbf{zeros} &   \mathbf{A[D(U)]}\\
	\end{array}
	\right).
$$
\end{defn}

\vspace{0.3cm}
\noindent \textbf{4.2 \textcolor{blue}{The Principal - natural  identifications}.} - see 4.1 in [6].

\vspace{0.1cm}
\noindent  Any \textbf{KoDAG}  is  a \textbf{di}-bicliques chain  $\Leftrightarrow$  Any  \textbf{KoDAG is a natural join} of complete bipartite \textbf{graphs} [ \textbf{di}-bicliques ]  = 
$$
	( \Phi_0 \cup \Phi_1 \cup ... \cup \Phi_n \cup ..., E_0 \cup E_1\cup ... \cup E_n \cup ...) \equiv D(\bigcup_{k\geq 0}\Phi_k,\bigcup_{k\geq 0} E_k )  \equiv D (\Phi,E) 
$$

\noindent where $E_k = \Phi_k\times \Phi_{k+1} \equiv \stackrel{\rightarrow}{K_{k,k+1}}$ and $E = \bigcup_{k\geq 0}E_k$.

\vspace{0.2cm}
\noindent Naturally, as indicated earlier in [6] any graded posets' Hasse diagram  with minimal elements set  including  \textbf{KoDAGs}  is of the form  
$$
	D (\Phi, E)  \equiv D( \bigcup_{k\geq 0}\Phi_k,\bigcup_{k\geq 0} E_k ) \Leftrightarrow  \langle \Phi,\leq \rangle
$$

\noindent where  $E_k \subseteq \Phi_k\times \Phi_{k+1} \equiv  \stackrel{\rightarrow}{K_{k,k+1}}$ and the  definition of  $\leq$  from  \textbf{1.3.} in [6] is applied.
\noindent In front of all the above presentation in [6] the following is clear .

\begin{observen}{(9 in [6].)} 
"Many"   graded digraphs with minimal elements set including \textbf{KoDAGs} $D = (V,\prec\!\!\cdot )$ denominate [encode] biunivoquely their  correspondent $(n+1)$-ary relation:  $\os_{k=0}^n E_k \subset \times_{k=0}^n \stackrel{\rightarrow}{K_{k,k+1}}$\\
where  $E_k \subseteq \Phi_k\times \Phi_{k+1} \equiv  \stackrel{\rightarrow}{K_{k,k+1}}$ and $n \in N \cup \{\infty\}$. 
\end{observen}
\vspace{0.2cm}
\noindent Which are those "many"? The characterization is arrived at   with au rebour point of view.
Any $n$-ary relation  ($n \in N \cup \{\infty\}$)  determines uniquely [may be identified with] its correspondent graded digraph with minimal elements set $\Phi_0$ given by the ($n$-ary rel.)  formula, where the sequence of binary relations $E_k \subseteq \Phi_k\times \Phi_{k+1} \equiv \stackrel{\rightarrow}{K_{k,k+1}}$ is denominated by the source $n$-ary relation as the following example shows.

\vspace{0.2cm}
\noindent \textbf{Example} (ternary = $Binary_1$ $\os$ $Binary_2$)

\noindent Let  $T \subset X\times Z\times Y$  where  $X =\{ x_1,x_2,x_3\}$, $Z = \{ z_1,z_2,z_3,z_4\}$, $Y = \{y_1,y_2\}$  and  
$$
	T = \{ \langle x_1,z_1,y_1 \rangle, \langle x_1,z_2,y_1 \rangle,   \langle x_1,z_4,y_2 \rangle, \langle x_2,z_3,y_2 \rangle, \langle x_3,z_3,y_2 \rangle \}.
$$

%%%%%%% picture 
\begin{figure}[ht]
\begin{center}
	\includegraphics[width=50mm]{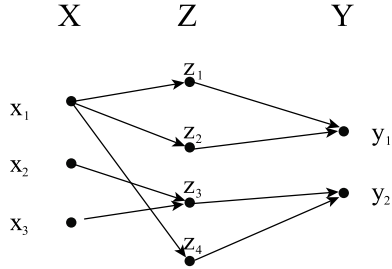}
	\caption {Display of example ternary = $Binary_1$ $\os$ $Binary_2$. \label{fig:ternary}}
\end{center}
\end{figure}

\vspace{0.1cm}

\begin{figure}[ht]
\begin{center}
	\includegraphics[width=100mm]{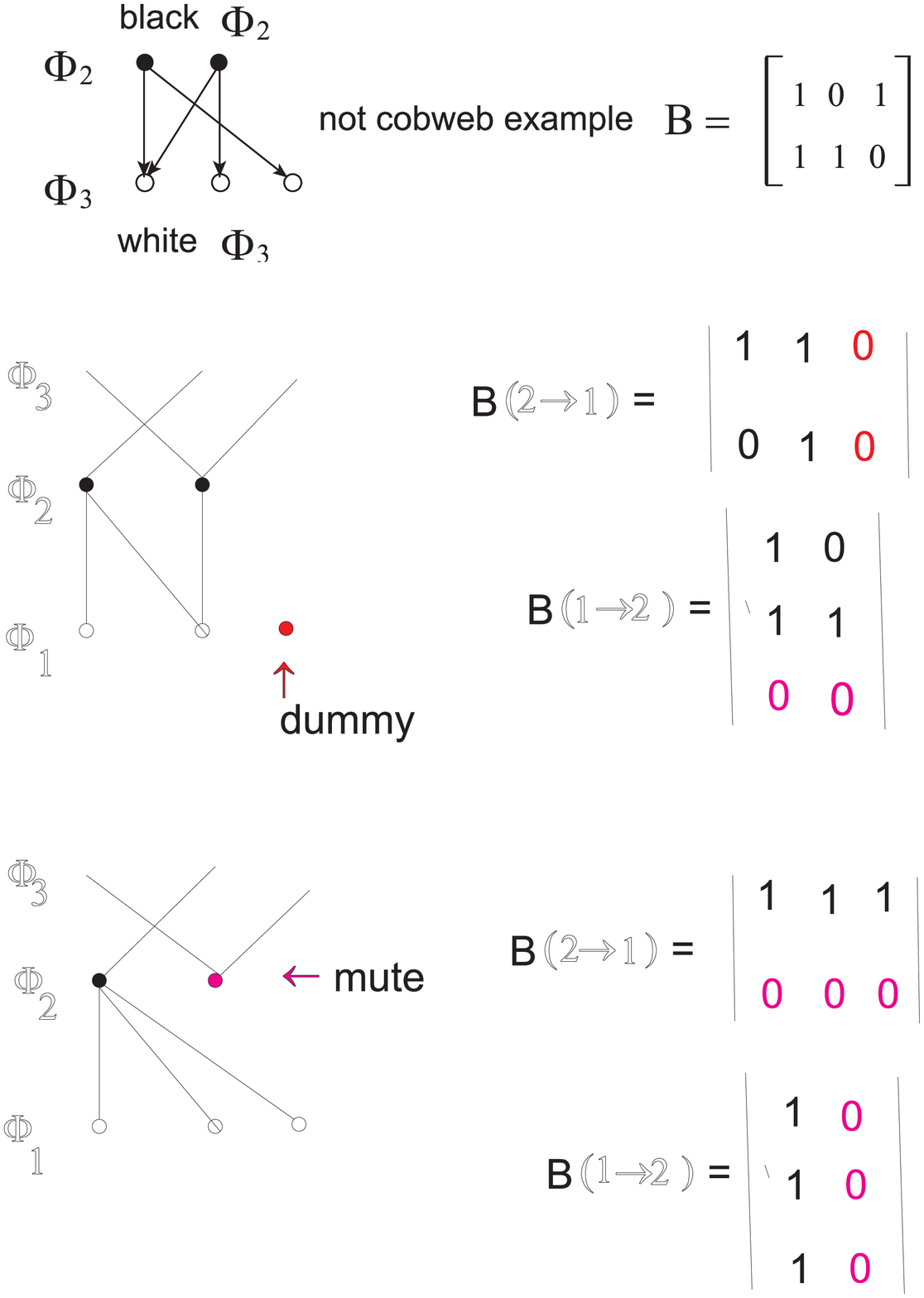}
	\caption {Display of bipartite digraphs with dummy and mute node and their biadjacency matrices. }\end{center}\label{fig:representation}
\end{figure}

\noindent Let $X\times Z \supset E_1= \{ \langle x_1,z_1 \rangle, \langle x_1,z_2 \rangle, \langle x_1,z_4 \rangle, \langle x_2,z_3 \rangle, \langle x_3,z_3 \rangle \}$  and 
$Z\times Y \supset E_2 = \{ \langle z_1,y_1 \rangle, \langle z_2,y_1\rangle,  \langle z_3,y_1\rangle, \langle z_4,y_2\rangle \}$.
Then $T = E_1 \os E_2$.\\

\vspace{0.1cm}
\noindent Finally - anyhow  we already know the now easy answer as it was expressed earlier.

\vspace{0.1cm}

\noindent $F$-\textit{graded poset may be identified with} $n$-\textit{ary relation  as above iff it is} $F$-\textit{graded poset with neither mute no dummy nodes.}\\ 
\vspace{0.1cm} 
\noindent Equivalently - \textit{zero columns or rows in bi-adjacency matrices of bipartite natural join summands of $P$ are forbidden}.  See and compare with relevant Fig. 2 at the beginning and Fig.17 at the end of the paper and perhaps consult [2,6].\\

\vspace{0.2cm}

\noindent \textbf{4.3 \textcolor{blue}{The Principal properties of cobweb posets as being also the natural join operation product}.}\\

\noindent It is well known  that any order relation is the intersection of its linear extensions as proved by Szpilrajn in  [60]. A set of linear extensions of $P$ whose intersection is $P$ is called a realizer of $P$. The dimension of an ordered set $P$ is the minimum cardinality of a realizer of $P$. \\
\noindent \textbf{The dimension of any cobweb poset is two}, as proved in [6]   (Observation 2  p.6).   
\noindent There it was also shown (Observation 3 p.7 in [6] )  that the \textbf{Ferrers dimension of any cobweb poset is one}.\\

\begin{lemma}{Some Properties of Cobwebs}\\ 
Cobweb posets are ranked, $N$-free $\Longleftrightarrow$ series-parallel posets.  Hence cobweb posets are greedy posets. Hence cobweb posets are \textbf{not} D-posets i.e. Dilworth posets.
\end{lemma}
 
\noindent Indeed. See it obvious in view what follows in more detail. \\

\noindent \textbf{4.1.} It is obvious that cobweb posets are series-parallel posets $\Longleftrightarrow$  $N$- free posets where:

\begin {defn} {[61]]}
A series-parallel poset is one that can be  recursively constructed by applying the operations of disjoint sum and linear product, starting with a single point .                 
\end{defn}
Note.  A series-parallel poset  $\Longleftrightarrow$   N- free poset which is is a particular case of  the $S$-poset  introduced in [62] by  Richard P. Stanley. See also p. 150 in [59].
\noindent (See - find  \textit{N poset's}  Hasse digraph in Fig.14).

\begin{lemma}( Characterization Lemma - see [63, 64, 65])
A finite poset is series-parallel if and only if it is N-free.  
\end{lemma}
We shall now mostly follow [65].

\vspace{0.2cm}

\noindent \textbf{4.2.}   From  Lemma 1.1 in [65] , namely : "A poset which does not have a subposet isomorphic to a cycle is a D-poset." we readily infer that cobweb posets are not D-posets i.e. Dilworth posets - (ad Dilworth posets  see also p.247 in [66]).

\vspace{0.2cm}

\noindent  \textbf{4.3.}  Cobweb posets are greedy  posets ( for greedy posets see: [65]  see also  p. 150 in [59] , see also p.242 in  [66] )
\begin {defn} {[65]]}
A linear extension $L = x_1, x_2,...,x_n$  of   $ P $ is greedy if  $L$ can be obtained by applying the following algorithm :  

\noindent 1. Choose a minimal element   $x_1 \in P$\\

\noindent 2. Suppose $x_1, ..., x_i$   have been chosen.  If there is at least one minimal element of $P \ \left\{ x_1, ..., x_i \right\} $    
which is greater than $x_i $ then  choose  $x_{i+1}$  to be any such minimal element; otherwise, choose $x_{i+1}$  to be any  minimal element of    $P \ \left\{ x_1, ..., x_i \right\}$.  
\end{defn}

\vspace{0.1cm}

\noindent  Let  $G(P)$   be the set of all greedy linear extensions of  $ P$ [65].   Let  $O(P)$ be the set of all optimal  linear extensions of  $ P$ [65]. ($L$  is called an optimal linear extension of $P$  if   $s(L, P) = s(P)$ where the jump number $s(P)$  of $P$  is defined as the minimum of  $s(L, P)$   over all linear extensions $L$ of  $P$. )

\vspace{0.1cm}

\noindent  Since the greedy algorithm above is a particular way of carrying out the algorithm for a linear extension, every poset $P$ has a greedy linear extension. 

\noindent Cogis [67] and Rival and Zaguia [68] had shown that every poset has an optimal greedy linear extension.

\vspace{0.2cm}

\noindent A poset P is greedy if G(P)   O(P) , that is, every greedy extension is optimal.

\vspace{0.2cm}

\noindent A poset P is N-free if P contains no cover -preserving subposet isomorphic to the poset N.

\vspace{0.1cm}

\noindent (See - find  \textit{ N poset's}  Hasse digraph in Fig.14).

\begin{lemma}( Lemma 1.2  in  [65])
Every N-free poset is greedy.
\end{lemma}

\vspace{0.1cm}

\noindent In 1982  El-Zahar and Rival [67] had proved that: 

\begin{lemma}( Lemma 1.3   in  [65])
A poset which does not contain a subposet isomorphic to Crown Poset satisfies O(P)   G(P).
\end{lemma}

\vspace{0.1cm}

\noindent Recall that the dual of the poset $P$ is the poset $P^d$  obtained from $P$ by reversing the order. 

\vspace{0.2cm}

\noindent Recall that a poset $P$  is said to be reversible if $L^d   G(P^d)$ for every $L   G(P)$.

\vspace{0.2cm}

\noindent\textbf{ Important.} In 1986  Rival and Zaguia [68] had proved the following characterization lemma.

\begin{lemma}( Lemma 1.4   in  [65])
A poset P is reversible if and only if O(P) = G(P).
\end{lemma}

\vspace{0.2cm}

\noindent \textbf{Exercise.}  Are cobweb posets reversible?

\vspace{0.2cm}

\noindent \textbf{Some open problems.}  For some open general problems see [5] on cobweb posets digraphs' elementary properties and questions.
For some open tiling task problems  see Maciej Dziemia\'nczuk's productions  [49] and [48].

\vspace{0.3cm}

\noindent \textbf{Acknowledgments}
\vspace{0.1cm}
\noindent Thanks are expressed here to the Student of Gda\'nsk University Maciej Dziemia\'nczuk for applying his skillful   TeX-nology with respect most of my articles since three years as well as for his general assistance and cooperation on KoDAGs  investigation. 

\vspace{0.1cm}

\noindent The author thanks  Dr Ewa Krot-Sieniawska for her critical remarks concerning mute nodes notion. The author expresses also his gratitude  also   Dr Ewa Krot-Sieniawska for her several years' cooperation and vivid application  of the alike  material deserving  Students' admiration for her being such a comprehensible and reliable  Teacher before she as Independent Person was fired by local Bialystok University local authorities exactly on the day she had defended  Rota and cobweb posets related dissertation with distinction.

\end{document}